\documentclass{elsart}
\journal{Stochastic Processes and their Applications}
\usepackage{epsfig, float}
\usepackage{graphicx}

\usepackage{amssymb,amstext,latexsym, amsmath, amscd, amsfonts, cite}
\usepackage[english]{babel}
\usepackage[latin1]{inputenc}

\usepackage{euscript}

\newtheorem{theorem}{Theorem}[section]
\newtheorem{propo}[theorem]{Proposition}

\newtheorem{lemma}[theorem]{Lemma}

\newtheorem{rmk}[theorem]{Remark}

\RequirePackage[colorlinks=true,citecolor=blue,linkcolor=blue]{hyperref}
\usepackage{hyperref}
\renewcommand{\theequation}{\thesection.\arabic{equation}}

\newcommand{\beq} {\begin{eqnarray*}}
\newcommand{\eeq} {\end{eqnarray*}}
\newcommand{\trm} {\textrm}
\newcommand{\tbf} {\textbf}
\newcommand{\tit} {\textit}
\newcommand{\noi} {\noindent}

\newcommand{\dis}{\displaystyle }

\def \P{\mathbb P}
\def \R{\mathbb R}
\def \N{\mathbb N}
\def \Z{\mathbb Z}
\def \E{\mathbb E}

\def\eps{\varepsilon}

\newcommand{\1}{{1\hspace{-0.2ex}\rule{0.12ex}{1.61ex}\hspace{0.5ex}}}

\begin{document}

\begin{frontmatter}

% Title, authors and addresses

% use the thanksref command within \title, \author or \address for footnotes;
% use the corauthref command within \author for corresponding author footnotes;
% use the ead command for the email address,
% and the form \ead[url] for the home page:
% \title{Title\thanksref{label1}}
% \thanks[label1]{}
% \author{Name\corauthref{cor1}\thanksref{label2}}
% \ead{email address}
% \ead[url]{home page}
% \thanks[label2]{}
% \corauth[cor1]{}
% \address{Address\thanksref{label3}}
% \thanks[label3]{}

\author[label1]{Claudie Chabriac}
\author[label1]{Agn\`es Lagnoux\corauthref{cor1}}
\corauth[cor1]{\textit{Corresponding author}. \tit{Phone}: +335.61.50.46.11\\
\tit{Email address}: lagnoux@univ-tlse2.fr\\
\tit{URL}: http://www.lsp.ups-tlse.fr/Fp/Lagnoux}
\author[label1]{Sabine Mercier}
\author[label2]{Pierre Vallois}
\address[label1]{Institut de
Math\'ematiques de Toulouse, UMR 5219, Universit\'e Toulouse 2,\\ 5 Allées Antonio Machado, 31058 Toulouse, France}
\address[label2]{Institut  Elie Cartan,
Universit\'e de Lorraine, CNRS UMR 7502, INRIA, BIGS, Campus Sciences, BP 70239, Vandoeuvre-lès-Nancy Cedex, 54506, France}

\title{Elements related to the largest complete excursion of a reflected BM stopped at a fixed time. Application to local score.}

\begin{abstract}
\noi We calculate the density function of $\big(U^*(t),\theta^*(t)\big)$,
where $U^*(t)$ is the maximum over $[0, g(t)]$ of  a reflected Brownian motion $U$, where  $g(t)$ stands for  the last zero of $U$ before $t$, $\theta^*(t)=f^*(t)-g^*(t)$, $f^*(t)$ is the hitting time of the level $U^*(t)$, and  $g^*(t)$ is the left-hand point of the interval straddling $f^*(t)$. We also calculate explicitly the marginal density functions of
$(U^*(t)$ and $\theta^*(t)$. Let $U^*_n$ and $\theta^*_n$ be the analog of  $U^*(t)$ and $\theta^*(t)$ respectively where the underlying process $(U_n)$ is the Lindley process, i.e. the difference between a centered real random walk and its minimum. We prove that
$\dis \Big(\frac{U^*_n}{\sqrt{n}},\; \frac{\theta^*_n}{n}\Big)$ converges weakly to
$\big(U^*(1),\theta^*(1)\big)$ as $n\rightarrow \infty$.
\end{abstract}

\begin{keyword}
Lindley process \sep local score \sep  Donsker invariance Theorem \sep reflected Brownian motion \sep inverse of the local time \sep Brownian excursions.\\
MSC: 60 F 17 \sep 60 G 17 \sep 60 G 40 \sep 60 G 44 \sep 60 G 50 \sep 60 G 52 \sep 60 J 55 \sep 60 J 65.
\end{keyword}
\end{frontmatter}

%%%%%%%%%%%%%%%%%%%%%%%
\section{Introduction}
%%%%%%%%%%%%%%%%%%%%%%%

\noi {\bf 1.1} The local score is a probabilistic tool which is often used by molecular biologists to study sequences of either amino-acids or nucleotides as DNA. In particular its statistical properties allow to determine the most significant segment in a given sequence, see for instance \cite{KAl90} and \cite{Wat95}. Any position $i$ in the sequence is allocated a random value $\epsilon_i$. For example, $\epsilon_i$  can measure either physical or chemical property of the $i$-th amino acid or nucleotide of the sequence. It can also  code the similarity between two components of two sequences. It is assumed that $(\epsilon_i)_{i\geq 1}$ is a sequence of independent and identically distributed random variables. Rather than considering $(\epsilon_i)_{i\geq 1}$, it is more usefull to deal with:
%Associated with a sequence $(\epsilon_i)_{i\geq 1}$  of independent and identically distributed random variables (r.v.'s),
%
\begin{equation}\label{defsn}
S_n=\epsilon_1+\cdots+\epsilon_n \ \mbox{ for } n\geq 1\ ; \quad S_0=0.
\end{equation}
Obviously, $(S_n)$ is the random walk starting at $0$, with independent increments $(\epsilon_i)_{i\geq 1}$. Let us introduce:
\begin{equation}\label{i1k}
  \underline{S}_n=\min\limits_{0\leq i\leq n}S_i,\quad n\geq 0.
\end{equation}
The two following processes $(U_n)$ and $(\overline{U}_n)$ play an important role in the study of biological sequences. The first one is called the Lindley process and is defined as:
\begin{equation}\label{defuk}
U_n=S_n-\underline{S}_n=S_n-\min\limits_{i\leq n}S_i, \quad n\geq 0.
\end{equation}
%
%$(U_n)$ is also known as the random walk $(S_n)$ reflected at its past minimum.\\\\
The process $(U_n)$ is non negative and further properties can be found either in (Chap. III of \cite{Asmussen03}) or Chap. I \cite{Borovkov76}).  The local score $\overline{U}_n$  is the supremum of the Lindley process up to time $n$.

\noi
 Molecular biologists are interested in "unexpected" large values of $(U_n)$, see \cite{Wat95}.

\noi The exact distribution of $\overline{U}_n$ has been determined in \cite{MD01}, using the exponentiation of a suitable matrix and classical tools related to Markov chains theory. Although the given formula in \cite{MD01}  is efficient whatever the sign  of $\E(\epsilon_i)$, in practice, it can be only applied to short sequences.
\noi However,  we are sometimes faced with long sequences and in these situations it is often assumed that they have  a negative trend, i.e. $\E(\epsilon_i)<0$. Then,  the local score $\overline{U}_n$ grows as $\ln(n)$ (see \cite{WGA87}) and an asymptotic approximation of the distribution of $\overline{U}_n$ as $n$ is large has been given in \cite{KAl90}, \cite{DK91},  using the renewal theory. When $\E(\epsilon_n)=0$, the asymptotic behavior of the tail distribution of $\overline{U}_n$ has been determined in \cite{DEV03} and the rate of convergence is given in \cite{EVa03}.

\noi Although the study of biological sequences is the starting point of this paper, the remainder will only consider the probabilistic model.

\noi
Here we consider that the $(\epsilon_i)_{i \geqslant 1}$ are centered with unit variance.

 \noi It is clear that the trajectory of  $(U_n)$ can be composed of a succession of $0$ and excursions above $0$.
%These excursions have an important biological interpretation. In particular, (CE QUI SUIT N'EST PAS TRES CLAIR) \textcolor{blue}{the highest one can corresponds to the best segment due to the physico chemical property or similarity scores that have been chosen. Note that the local score  $\overline{U}_n$ can be viewed  as the maximum of the heights of all the excursions up to time $n$.}
However,  we only deal with  \tit{complete} excursions up to a fixed time. This leads us to introduce
the maximum $U^*_n$ of the heights of all the complete excursions up to time $n$. The second variable which will play an important role is $\theta^*_n$, the time necessary to reach  its maximum height $U^*_n$. See Section \ref{sec:approx_discret} for more informations and detailed definitions of the previous RVs.

\noi We believe that the knowledge of  the joint distribution of the pair $(U^*_n,\theta^*_n)$ should permit the associated bi-dimensional statistical tests to be more powerful  than the usual ones based on the first component. This program  should be developed in a forthcoming paper.

\noi
%The motivation of the work proposed in this article is to take into account the length of the segment that realized the local %score. In pratice, we determine the joint distribution of $U^*_n$ and $\theta^*_n$ as a first step. \\

\noi {\bf 1.2}\  Unfortunately, it is difficult to determine explicitly the law of $(U^*_n,\; \theta^*_n)$ for a fixed $n$. Bearing in mind applications with  long   biological sequences, it is relevant to study the distribution of $(U^*_n,\; \theta^*_n)$ where $n$ is large. The functional convergence theorem  of Donsker tells us  that the initial random walk $(S_k,\; 0\leq k\leq n)$ normalized by the factor $1/\sqrt{n}$ converges in distribution as $n\rightarrow\infty$ to the Brownian motion $(B(s),\; 0\leq s \leq 1)$, see  Sections \ref{sec:link} and \ref{sec:cv_proc} for a more precise formulation. It is easy to deduce that the normalized Lindley process $\dis \Big(\frac{U_k}{\sqrt{n}}, \; 0\leq k\leq n\Big)$ can be approximated by  $(\widehat{U}_s,\; 0\leq s\leq 1)$ where:
 %the so-called reflected Brownian motion
%
\begin{equation}\label{mbr7}
    \widehat{U}(t):=B(t)-\inf_{0\leq s\leq t}B(s),\quad t\geq 0.
\end{equation}
\noi
Recall that the process $(\widehat{U}(s),\; s\geq 0)$ is distributed as the reflected Brownian motion, since:
 \begin{equation}\label{trL1}
   \big(|B(t)|, \;t\geq 0\big) \ \overset{(d)}{= } \ \big(B(t)-\min_{0\leq u\leq t}B(u),\; t\geq 0\big).
\end{equation}
It turns out that the asymptotic behavior of $(U^*_n,\; \theta^*_n)$ for large $n$ should be closely linked the distribution of $\big( U^*(1), \theta^*(1)\big)$ where $ U^*(1)$ and $\theta^*(1)$ are the analog in continuous time of
   $U^*_n$ and $\theta^*$. Consequently, the knowledge of the distribution of $(U^*_n,\; \theta^*_n)$ for large $n$ reduces to $\big(U^*(1),\theta^*(1)\big)$. Let us briefly define these RVs.  As we proceed with the random walk $(S_n)$, we introduce the following processes (see Section \ref{sec:thal_results} for more explicit definitions):

 %leads us to study the analog of
 %$(U^*_n,\; \theta^*_n)$ in the setting of continuous time process $(U(t))$. We proceed similarly as for the Lindley process and define

\begin{enumerate}
\item the local score $\overline{U}(t)$ which is the maximum of the heights of all the excursions of $U(s)$ up to time $t$, i.e. $\dis U^*(t):=\sup_{0\leq s\leq t}U(s)$,
\item  the maximum $U^*(t)$ of the heights of all the \tit{complete} excursions up to time $t$,
\item    the time $\theta^*(t)$ taken by $U$ to reach $U^*(t)$ starting from the beginning of this highest excursion.
\end{enumerate}

%\noi  The  approximation of $(U_n,\; n\geq 0)$ by $(\widehat{U}(t),\; t\geq 0)$ allows to prove that the asymptotic %distribution of $\dis \Big(\frac{U^*_n}{\sqrt{n}},\; \frac{\theta^*_n}{n}\Big)$ as $n\rightarrow \infty$  is the one of %$\big(U^*(1),\theta^*(1)\big)$.

\noi

\noi {\bf 1.3} Let $t$ be a fixed real number. The density function of $\overline{U}(t)$ is known (see either Subsection 2.11 in  \cite{Billingsley99} or Lemma 3.2 in \cite{RVY08}). Although $U^*(t)=\overline{U}\big(g(t)\big)$ and $g(t)$ is not a stopping time it is however easy to calculate  the density function  of $U^*(t)$. Indeed, the process $\left(g(t)^{-1/2}B(g(t)s),\; 0\leq s\leq 1\right)$ is distributed as $\big(b(s),\; 0\leq s\leq 1\big)$ and is independent of $g(t)$, where $g(t)$ is the last zero of $U(s)$ before $t$ and $b$ is the Brownian bridge (see e.g. \cite{Bertoin96}). Therefore:
\begin{equation}\label{16h24}
    U^*(t)\overset{(d)}{=} \sqrt{g(t)}\sup_{0\leq s\leq 1} |b(s)|.
\end{equation}
Finally, we conclude using the fact that $g(t)$ is distributed as the arcsine law (see again \cite{Bertoin96}) and the distribution of $\sup_{0\leq s\leq 1} |b(s)|$ is given by the Kolmogorov-Smirnov formula   (see e.g. \cite{PY01}). The final and explicit result is given in Theorem \ref{th:marginals}.

\par
\noi However, as far as we know, the distribution of $\big(U^*(t),\theta^*(t)\big)$ is unknown.
\noi Using the  theory of excursions related to the one dimensional Brownian motion, we determine in Theorem \ref{th:th1} the density function of the couple $\big(U^*(t),\theta^*(t)\big)$. Since  we are interested in statistical tests based on the joint law of $(U^*(t),\theta^*(t))$, then   we have to determine the quantiles of $(U^*(t),\theta^*(t))$. Unfortunately the expression of the density function is complicated and does not allow us to calculate the distribution function of $(U^*(t),\theta^*(t))$. In Theorem \ref{th:th2}, we express, for any bounded Borel function $f :\; ]0,\infty[\times ]0,\infty[\rightarrow \mathbb{R}$, the expectation of    $f\big(U^*(t),\theta^*(t)\big)$ as $\E(f(A_1)A_2)$ where $A_1$ and $A_2$ are RVs which can be simulated. Therefore, the quantity  $\E\big( f\big(U^*(t),\theta^*(t)\big)$ can be approximated by a Monte-Carlo scheme.

 %are r.v.s with known distribution. and the probability $P\big(U^*(t)\leq a ,\theta^*(t)\leq b\big)$  for any $a, b>0$ can be computed via Monte Carlo simulations.\par

\noi  In Section \ref{ssec:pair} we fix notations related to the setting of processes in continuous time, i.e. here the underlying process is the Brownian motion. The main results are Theorems \ref{th:th1}, \ref{th:th2}, \ref{th:marginal_theta} and \ref{th:marginals} and they  are given in Section \ref{ssec:marges}.  Theorem \ref{th:th1} is based on Propositions \ref{prop:prop2} and \ref{prop:prop1}. Although the law of $U^*(t)$ is easy to calculate, $\theta^*(t)$ is more difficult, see Theorem \ref{th:marginal_theta}.  We recall in Section \ref{sec:approx_discret} the functional  approximation of the one dimensional Brownian motion by normalized random walks. Then, with additional technical developments, see Proposition \ref{prop:link} and Theorem \ref{T1} we  obtain the weak convergence of
$\dis \Big(\frac{U^*_n}{\sqrt{n}},\; \frac{\theta^*_n}{n}\Big)$ as $n\rightarrow \infty$  towards  $\big(U^*(1),\theta^*(1)\big)$. All the proofs which are not immediate have been given in Section \ref{Pro1}.\\

%\ref{sec:thal_results}, we start giving theoretical results, i.e. whose related to Brownian motion and in particular the distribution of the pair $\big(U^*(t),\theta^*(t)\big)$. The marginal laws of $U^*(t)$ and $\theta^*(t)$ are also made explicit and the one of $U^*(t)$ is identified with the one given above.

\noi {\bf Acknowledgements} The authors are greatly indebted to the referee for his fruitful comments, references and suggestions.

%%%%%%%%%%%%%%%%%%%%%%%%%%%%%%%%%%%%%%%%%%%%%%%%%%%%%%%%%%%%%%%%%%%%%%%%%%%%%%%%%%%%%%%%%%%%%%%%%%%%%%%%%%%%%%%%%
%%%%%%%%%%%%%%%%%%%%%%%%%%%%%%%%%%%%%%%%%%%%%%%%%%%%%%%%%%%%%%%%%%%%%%%%%%%%%%%%%%%%%%%%%%%%%%%%%%%%%%%%%%%%%%%%%
%%%%%%%%%%%%%%%%%%%%%%%%%%%%%%%%%%%%%%%%%%%%%%%%%%%%%%%%%%%%%%%%%%%%%%%%%%%%%%%%%%%%%%%%%%%%%%%%%%%%%%%%%%%%%%%%%
%%%%%%%%%%%%%%%%%%%%%%%%%%%%%%%%%%%%%%%%%%%%%%%%%%%%%%%%%%%%%%%%%%%%%%%%%%%%%%%%%%%%%%%%%%%%%%%%%%%%%%%%%%%%%%%%%

%%%%%%%%%%%%%%%%%%%%%%%
\section{Theoretical results}\label{sec:thal_results}
%%%%%%%%%%%%%%%%%%%%%%%

\setcounter{equation}{0}

\subsection{Notation}\label{ssec:nota}

\noi Let $(B(t),\; t\geq 0)$ be a standard Brownian motion started at 0 and $U(t)$ is the reflected Brownian motion at time $t$:
\begin{equation}\label{tr1}
    U(t):=\vert B(t)\vert, \quad t\geq 0.
\end{equation}

 \noi
The excursion (above $0$) straddling $t$ starts at $g(t)$ and ends at $d(t)$, namely
 \begin{equation}\label{tr2}
    g(t)=\sup\{s\leq t,\ U(s)=0\},\quad d(t)=\inf\{s\geq t,\ U(s)=0\},\quad t\geq 0.
 \end{equation}
%
%%is the analog in continuous time of the Lindley process
%\begin{equation}\label{mbr7}
 %   U(t):=B(t)-\inf_{0\leq s\leq t}B(s),\quad t\geq 0.
%\end{equation}
%
 \noi Let $\overline{U}(t)$ be the supremum of $U$ over $[0,t]$
\begin{equation}\label{deHc1}
    \overline{U}(t):=\sup_{0\leq s\leq t}U(s),\quad t\geq 0.
\end{equation}
\noi
%Since $U$ is a continuous process, there exists a random time $t^*$ in $[0,t]$ such that %$\overline{U}(t)=\overline{U}(t^*)$. Let $g(t)=\sup\{s\leq t,\ U(s)=0\}$ and $d(t)=\inf\{s\geq g(t),\ U(s)=0\}$.
% If $d(t)\leq t$, this excursion is achieved before $t$ and is complete. Otherwise, the highest complete excursion is %the biggest one before $g(t)$. Finally
Then, the highest height $U^*(t)$ of all the {\it complete} excursions of the process $(U(r)\ ;\ 0\leq r\leq t)$ equals
%up to time $t$ and $\theta^*(t)$ the associated length are defined as follows:
%
\begin{equation}\label{formule_U}
U^*(t):=\overline{U}(g(t))=\sup\limits_{0\leq s\leq g(t)}U(s),\quad t\geq 0.
\end{equation}
\noi Let $f^*(t)$ be the unique time which achieves the maximum of $U$ over $[0,\; g(t)]$:
\begin{equation}\label{tr3}
    f^*(t):=\sup\{r\leq g(t)\ ;\ U(r)=U^*(t)\},\quad t\geq 0.
\end{equation}
It is worth introducing  the left end-point $g^*(t)$ of the  excursion straddling $ f^*(t)$:
\begin{equation}\label{formule_gf}
g^*(t):=g(f^*(t))=\sup\{r\leq f^*(t)\ ;\ U(r)=0\},\quad t\geq 0
\end{equation}
as well as   the right end-point $d^*(t)$ of this  excursion:
\begin{equation}\label{formule_d}
d^*(t):=d(f^*(t))=\inf\{r\geq f^*(t)\ ;\ U(r)=0\},\quad t\geq 0.
\end{equation}
\noi It is convenient to visualize the different variables in Figure 1. %\ref{LabelNotationsUt}.\\

\begin{figure}[h]\label{LabelNotationsUt}
\centering
\includegraphics[scale=0.3]{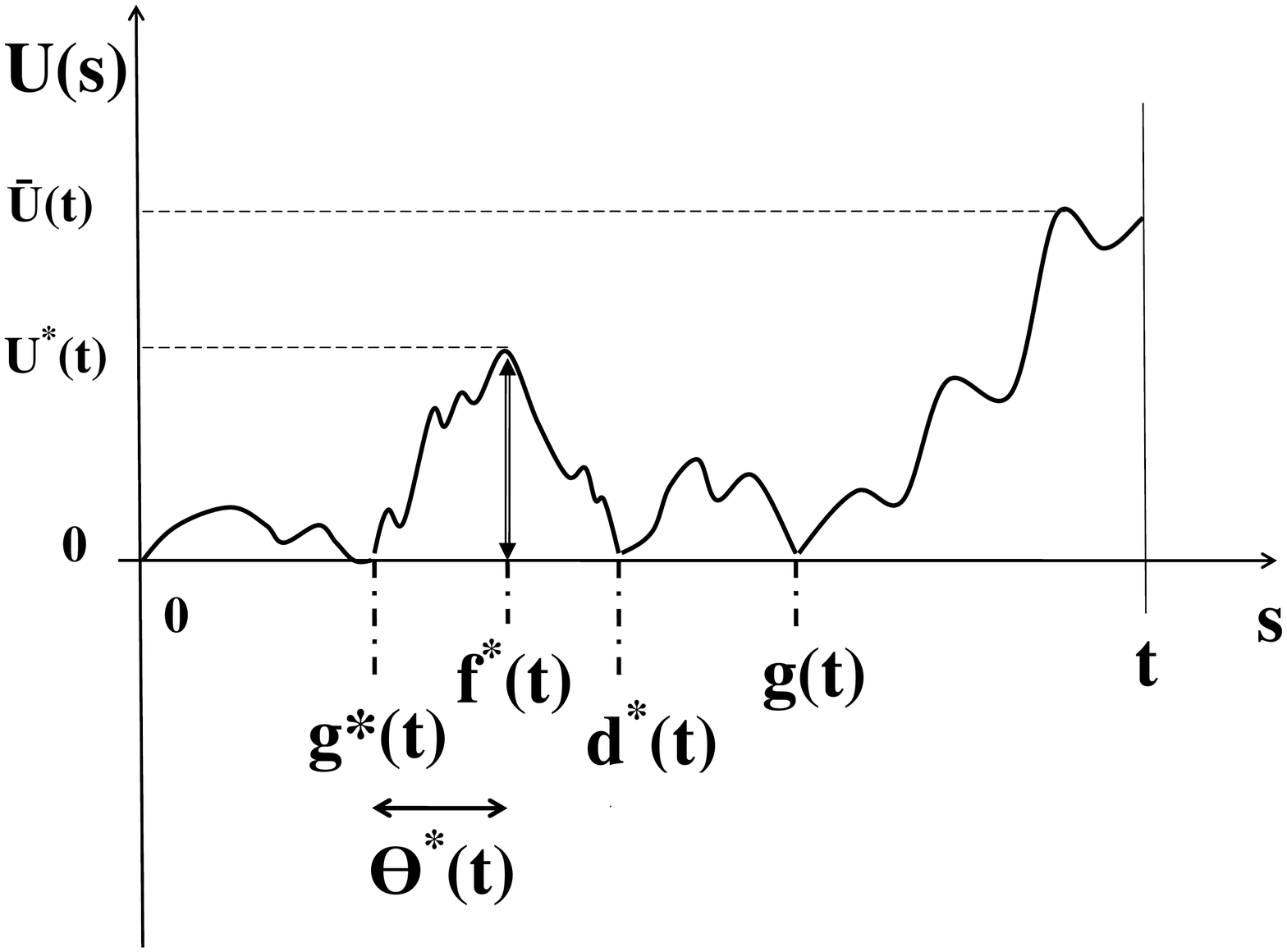}
\caption{Notation $U(t)$}
\end{figure}

\noi
We are interested in the joint law of $U^*(t)$ and $\theta^*(t)$ where the second variable is defined as
\begin{equation}\label{tr8}
\theta^*(t):=f^*(t)-g^*(t), \quad t\geq 0.
\end{equation}

It is convenient to introduce the following notation which will be used extensively in the sequel.

%%%%%%%%%%%%%%%
%\begin{nota}\label{ntr1}
%
1) $(\xi_n)_{n\geq 1}\cup\{\xi,\xi'\}$ is a family of i.i.d. r.v.s such that

\begin{equation}\label{formule_xi}
\xi\overset{(d)}{=}\xi'\overset{(d)}{=}\xi_n\overset{(d)}{=}T_1(R)
\end{equation}
with
\begin{equation}\label{htBe1}
 T_x(R)=\inf\{s\geq 0\ ;\ R(s)=x\}, \quad x>0
\end{equation}
and $(R(s), s\geq 0)$ stands for a  3-dimensional Bessel process started at 0.\\
%
%We assume that $\xi$, $\xi'$ and $U$ are independent.\\\\
%
The density $p_\xi$ is explicitly known and is given by

\begin{eqnarray}
p_{\xi}(u)&=&\frac{1}{\sqrt{2\pi}u^{3/2}}\sum_{k\in\Z}\left (-1+\frac{(1+2k)^2}{u}\right ) \exp{\left ( -\frac{(1+2k)^2}{2u}\right )}\label{eq:dens_xi}\\
&=& \frac{d}{du}\left (\sum_{k\in\Z} (-1)^k \exp\left (-\frac{k^2\pi^2u}{2}\right )\right )\label{eq:dens_xi_2}
\end{eqnarray}

%\cite{CieTay} or formula 2.0.2 in \cite{BS02}). 
(see for instance \cite{BPY01} p 8 and 24). In perspective of simulation, let us mention that an efficient algorithm to simulate very quickly the r.v. $\xi$ is given in \cite{DH11}.

2) $e_0'$, $(e_n)_{n\geq 0}$ is a sequence of i.i.d. exponential r.v.s.

%We assume moreover that $(e_n)_{n\geq 1},$ $(\xi_n)_{n\geq 1},$ $\xi$ and $\xi'$ are independent.

3) $(\lambda(x), x\geq 0)$ is the  process defined by
\begin{equation}\label{formule_lambda}
\lambda(x):=x^2(\xi_1+\xi_2)+\sum\limits_{k\geq 1}\frac{\xi_{2k+1}+\xi_{2k+2}}{\left(\frac 1 x+e_1+\cdots+e_k\right)^2},\quad x\geq 0.
\end{equation}
The sum converges a.s. and in $L^1$ (see Lemma \ref{lem:lem_cv_lambda}). The Laplace transform of $\lambda(x)$ has been calculated in \eqref{eq:laplace}.

4) $\alpha_1$ and $\alpha_2$ are two $[0,1]$ valued r.v.s; $\alpha_2$ is uniformly distributed and the density function of $\alpha_1$ is $\frac{2}{\pi}\frac{1}{\sqrt{1-s^2}}\1_{[0,1]}(s)$.\\\\
We always assume in the sequel that 
\begin{equation}\label{tr9}
   e_0', \; (e_n)_{n\geq 0}, \; (\xi_n)_{n\geq 1},\; \xi,\; \xi',\; \alpha_1, \; \alpha_2 \mbox{ and } (U(t))_{t\geq 0} \mbox{ are independent.}
\end{equation}

%\end{nota}
%%%%%%%%%%------------

\subsection{Distribution of the pair $(U^*(t),\theta^*(t))$}\label{ssec:pair}

\noi
%First we determine the density function of $(U^*(t),\theta^*(t))$ (Theorem \ref{th:th1}). 
%We also give a formula for the expectation of $f(U^*(t),\theta^*(t))$ which permits its estimation via a Monte Carlo algorithm (Theorem \ref{th:th2}).
%%\noi To give  either the density function or the expectation of a given function of $\big(U^*(t),\theta^*(t)\big)$ 
%In that view, it is convenient to introduce preliminary notation.

The main results are Theorems \ref{th:th1}, \ref{th:th2}, \ref{th:marginal_theta} and \ref{th:marginals}. All the proof of results stated in this section will be developed in Section \ref{Pro1}.\\
In Theorem \ref{th:th1}, we determine the density function of $(U^*(t),\theta^*(t))$. Its proofs is based on the theory of excursion related to the Brownian motion, see for instance Chap XII in \cite{RY99}. 
Let us briefly recall the ingredients which are needed. Let $(L(t),\; t\geq 0)$ be the local time process at $0$ related to the Brownian motion $(B(t),\; t\geq 0)$. The random function $t \mapsto L(t)$ is continuous and non-decreasing. Let $(\tau_s,\, s\geq 0)$ be its right inverse. The proof of Theorem \ref{th:th1} has two main steps. In Proposition \ref{prop:prop1} below we begin with expressing the distribution of $\big(U^*(t),\theta^*(t)\big)$ in terms of the one of $\big(\overline{U}(\tau_1),\; \tau_1\big)$.

%%%--------------------------------------------
%\begin{propo}\label{prop:prop1}
%Let $t$ be a fixed positive real number. Then, for any bounded Borel function  $f :\; \R_+^2\rightarrow \mathbb{R} $,
%%
%$$\mathbb{E}\left[f\big(U^*(t),\theta^*(t)\big)\right]=\sqrt{\frac{2}{\pi}}\int_0^{+\infty}\mathbb{E}\left[f(y,y^2\xi)\psi(y)\right]\,\frac{dy}{y^2}$$
%%
%where
%%
%\begin{equation}\label{formule_psi}
%\psi(y):=\frac{1}{\tau_1}\left[\sqrt{\left(t-y^2(\xi+\xi')\right)_+}-\sqrt{\left(t-y^2\left(\xi+\xi'+\frac{\tau_1}{(\overline{U}(\tau_1))^2}\right)\right)_+}\right],
%\end{equation}
%$x_+:=\sup\{x,0\}$.
%\end{propo}
%%%%%%%%%%%%%%%%%%
%%see the proof in Appendix \ref{app:th1} for more details).

%%--------------------------------------------
\begin{propo}\label{prop:prop1}
Let $t$ be a fixed positive real number. Then, the density function 
 of $\left(U^*(t),\theta^*(t)\right)$ is given by
\begin{equation}\label{eq:dens_couple}
p_{(U^*(t),\theta^*(t))}=\sqrt{\frac{2}{\pi}}\rho(x,y)p_{\xi}\left(\frac{y}{x^2}\right)\frac{1}{x^4}, \quad x>0,0<y<t
\end{equation}
where $p_{\xi}$ is the density function of $\xi$ (see \eqref{eq:dens_xi}-\eqref{eq:dens_xi_2}),
\begin{equation}\label{formule_rho}
\rho(x,y):=\int_0^{+\infty} \E\left(\frac{1}{\tau_1}\left\{\sqrt{\left(\rho_1\right)_+}-\sqrt{\left(\rho_1-\frac{x^2\tau_1}{(\overline{U}(\tau_1)^2}\right)_+}\right\}\right)p_{\xi}(u)du,
\end{equation}
$x_+:=\sup\{x,0\}$ and $\rho_1:=t-y-x^2u$.
\end{propo}
%%%%%%%%%%%%%%%%%

\noi
We are then naturally lead to determine the distribution  $\big(\overline{U}(\tau_1),\tau_1\big)$.

%%%-----------------------------------------------------
\begin{propo}\label{prop:prop2}
\begin{enumerate}
  \item For any $x\geq 0$, the sum in (\ref{formule_lambda})  converges a.s. and in $L^1$.
  \item The r.v. $\overline{U}(\tau_1)^{-1}$ is exponentially distributed and conditionally on\\ $\left\{\overline{U}(\tau_1)=x\right\}$, $x>0$,
\begin{equation}\label{formule_loi_tau_1}
\tau_1\overset{(d)}{=}\lambda(x).
\end{equation}
%
%where $\lambda(x)$ has been defined by \eqref{formule_lambda}.
\end{enumerate}
\end{propo}
%%-----------------------------------------------------------------

\noi Finally, combining Propositions \ref{prop:prop1} and \ref{prop:prop2} provides the density function of  $\big(U^*(t),\theta^*(t)\big)$.

%%%-----------------------------------
%\begin{theorem}\label{th:th1}
%For any $t>0$, the pair $\big(U^*(t),\theta^*(t)\big)$ has the following density
%\begin{equation}\label{formule_dens}
%p_{(U^*(t),\theta^*(t))}(x,y)=\sqrt{\frac{2}{\pi}}\frac{1}{x^4}\,p_\xi\left(\frac{y}{x^2}\right)\,\psi_1\left(x,\frac{y}{x^2}\right)
%\,\1_{\{0<x,0<y<t\}}
%\end{equation}
%%
%with
%%
%%\begin{eqnarray}
%%&&\psi_1(x,u):=\label{15j13}\\
%%&&\int_0^{+\infty}\mathbb{E}\left[\frac{1}{\lambda(v)}\left\lbrace\sqrt{\left(t-x^2(u+\xi')\right)_+}
%%-\sqrt{\left(t-x^2\left[u+\xi'+\frac{\lambda(v)}{v^2}\right]\right)_+}\right\rbrace \right] \frac{e^{-1/v}}{v^2}\,dv\nonumber
%%\end{eqnarray}
%\begin{eqnarray}
%&&\psi_1(x,u):=\int_0^{+\infty} \frac{e^{-1/v}}{v^2} \label{15j13}\\
%&&\times \mathbb{E}\left[\frac{1}{\lambda(v)}\left\lbrace\sqrt{\left(t-x^2(u+\xi')\right)_+}
%-\sqrt{\left(t-x^2\left[u+\xi'+\frac{\lambda(v)}{v^2}\right]\right)_+}\right\rbrace \right] \,dv\nonumber
%\end{eqnarray}
%%
%where $x>0$ and $0<u x^2<t$.
%\end{theorem}
%%%%------------------------------------------------------

%%-----------------------------------
\begin{theorem}\label{th:th1}
For any $t>0$, the pair $\big(U^*(t),\theta^*(t)\big)$ has the density function $p_{(U^*(t),\theta^*(t))}$ given by \eqref{eq:dens_couple} where
\begin{equation}\label{formule_rho_2}
\rho(x,y)=\int_{\R^2_+} \E\left(\frac {1}{\lambda(v)} \left\{\sqrt{\left(\rho_1\right)_+}-\sqrt{\left(\rho_1-\frac{x^2}{v^2}\lambda(v)\right)_+}\right\}\right)p_{\xi}(u)\frac{e^{-1/v}}{v^2}\,dudv
\end{equation}
and $\rho_1$ has been defined in Proposition \ref{prop:prop1}.
\end{theorem}
%%%------------------------------------------------------

\noi Formula \eqref{formule_rho_2} has the disadvantage to be not completely explicit and therefore it does not allow a direct calculation of $\mathbb{E}\left[f\big(U^*(t),\theta^*(t)\big)\right]$ for a given bounded Borel function $f$. For instance, for our biological motivation explained in the Introduction, it would be interesting to calculate $\mathbb{P}\big(U^*(t)\leq a,\ \theta^*(t)\leq b\big)$ for any $a, b>0$. Rewriting the proof of Theorem \ref{th:th1} leads us to an equivalent formulation of Theorem \ref{th:th1} which gives rise to a more useful formula.

%%----------------------------------
\begin{theorem}\label{th:th2}
%Let $t>0$. For any bounded Borel function $f: \; [0,\infty[\times [0,\infty[\rightarrow\mathbb{R}$ we have
%$$\mathbb{E}\left[f(U^*(t),\theta^*(t))\right]=\mathbb{E}\left[f\left(\frac{\sqrt{t}|G|}{e_0},\frac{tG^2\xi}{e_0^2}\right)\frac{e^{ \frac{G^2}{2}+e_0}}{\sqrt{1-G^2 \left(\frac{\xi+\xi'}{e_0^2}+\lambda\left(\frac{1}{e'_0}\right)\right)}}
%\right.\qquad \qquad \qquad \qquad \qquad \qquad \qquad \qquad  $$
%%
%\begin{equation}\label{15jui}
   %\qquad \qquad \qquad \qquad \qquad \qquad \qquad \qquad \qquad \left.\times \1_{\left\lbrace G^2 \left(\frac{\xi+\xi'}{e_0^2}+\lambda\left(\frac{1}{e'_0}\right)\right)\leq 1,\ e_0<e'_0\right\rbrace}\right]
%\end{equation}
%%
 %where the r.v.'s $e_0,\; e'_0,\; G,\; \xi$ and $\xi'$ are independent and their distributions are given by  \eqref{tr4}, \eqref{tr5} and \eqref{formule_xi}.
Let $f:\mathbb{R}_+^2\to\mathbb{R}$ be a bounded Borel function. Then
\begin{equation}\label{15jui}
\E\left(f(U^*(t),\theta^*(t))\right)=\sqrt{\frac{\pi}{2}}\E\left[f\left(\frac{\alpha_1\sqrt t}{\sqrt Z},\frac{t \alpha_1^2 \xi}{Z}\right)\frac{\alpha_2 e_0^{'2}}{\sqrt Z}\right]
\end{equation}

where $Z=\xi+\xi'+e_0^{'2}\alpha_2^2\lambda(1/e_0').$
\end{theorem}
%%------------------------

%\noi
%The proofs of Theorems \ref{th:th1} and \ref{th:th2} are postponed in Sections \ref{app:th1} and \ref{app:th2}. Theorem \ref{th:th2} comes directly by rewriting the proof of Theorem \ref{th:th1}.

\subsection{Distributions of $U^*(t)$ and $\theta^*(t)$}\label{ssec:marges}

We begin with the distribution of $\theta^*(t)$.

\begin{theorem}\label{th:marginal_theta}
For any $t>0$, $\theta^*(t)$ admits the following density function
\begin{equation}\label{eq:marginal_theta}
f_{\theta^*(t)}(x)=\frac{1}{x}\sum_{k\geq 1}(-1)^{k+1}\frac{\sinh\left (\pi k\sqrt{\frac{x}{t-x}}\right )}{\cosh^2\left( \pi k\sqrt{\frac{x}{t-x}}\right )} \1_{[0,t]}(x).
\end{equation}
\end{theorem}

%The proof of Theorem \ref{th:marginal_theta} is postponed in Section \ref{sec:th:marginal_theta}. 
We now consider the law of $U^*(t)$.

\begin{theorem}\label{th:marginals}
Let $t>0$.
%\begin{enumerate}
%\item Let $b(t)$ be the Brownian bridge at time $t$ and $g(t)$ the last zero of the process $U$ before time $t$. Then
%\begin{equation}\label{eq:dist_U*_PY}
%U^*(t)\overset{(d)}{=} \sqrt{t}\sqrt{g(1)}b^*,
%\end{equation}
%where $g_1$ and $b^*$ are two independent r.v.s such that
%
%\begin{eqnarray}
%\P(g(1)\in dx)&=&\frac{1}{\pi\sqrt{x(1-x)}}\1_{[0,1]}(x)dx \label{eq:dist_g1};\\
%\P(b^*>x)&=&2\sum_{k\geq 1} (-1)^{k-1} e^{-2k^2x^2}, \quad (x>0). \label{eq:dist_b*}
%\end{eqnarray}
\begin{enumerate}
\item We have the following identity in law
\begin{equation}\label{eq:dist_U*_PY}
U^*(t)\overset{(d)}{=} \sqrt{t}\sqrt{g(1)}b^*,
\end{equation}
where $g(1)$ and $b^*$ are two independent r.v.s such that

\begin{eqnarray}
\P(g(1)\in dx)&=&\frac{1}{\pi\sqrt{x(1-x)}}\1_{[0,1]}(x)\,dx \label{eq:dist_g1},\\
\P(b^*>x)&=&2\sum_{k\geq 1} (-1)^{k-1} e^{-2k^2x^2}, \quad x>0. \label{eq:dist_b*}
\end{eqnarray}

\item $U^*(t)$ admits the following density function

\begin{equation}\label{eq:def_phi}
f_{U^*(t)}(x)=4\sqrt{\frac{2}{\pi t}} \left(\sum_{k\geq 1} (-1)^{k-1} k e^{-\frac{2k^2x^2}{t}}\right), \quad x>0.
\end{equation}
\end{enumerate}
\end{theorem}

The proof of item 1 in Theorem \ref{th:marginals} is straightforward and has been developed in the Introduction. Note that this 
direct approach does not use the knowledge of the density function of 
$(U^*(t),\theta^*(t))$. However, Lemma \ref{lem:lem_f1} permits to get another expression of the distribution of $U^*(t)$.

%The proof of Theorem \ref{th:marginals} (see Section \ref{sec:th:marginals}) does not use our results related to the joint distribution of the pair $(U^*(t),\theta^*(t))$. However using Lemma \ref{lem:lem_f1} and the following result, one can compute in a different manner the density of $U^*(t)$ given by \eqref{eq:def_phi}.

\begin{propo}\label{prop:marginals}
We define two other stopping times $T_U(a)=\inf\{t\geq 0,\ U(t)=a\}$ and $T_{\widehat{B}}(a)=\inf\{t\geq 0,\ \widehat{B}(t)=a\}$ where $(\widehat{B}(t),\; t\geq 0)$ is a one dimensional Brownian motion independent of $U$. Then
\begin{equation}\label{eq:inverse_2}
\P(U^*(t)>a)=\P\left(T_U(a)+T_{\widehat{B}}(a)<t\right), \quad t>0, a>0.
\end{equation}
\end{propo}

\begin{rmk}
It is clear that \eqref{eq:inverse_2} allows to compute the cumulative distribution function of $U^*(t)$ and gives a complement to \eqref{eq:def_phi}. The distributions of  $T_{\widehat{B}}(a)$ and $T_{U}(a)$ are explicitly known: the density function of $T_{\widehat{B}}(a)$ is 
$\frac{a}{\sqrt{2\pi t^3}}e^{-a^2/2t} \1_{\{t>0\}}$
and (see Lemma 3.2 in \cite{V91})
$$\P\left(T_{U}(a)>t\right)=\frac{4}{\pi}\sum_{k\geq 0} \frac{(-1)^k}{2k+1}\exp\left\{-\frac{(2k+1)^2\pi^2}{8a^2}t\right\},\ t>0.$$
\end{rmk}

%%%%%%%%%%%%%%%%%%%%%%%%%%%%%%%%%%%%%%%%%%%%%%%%%%%%%%%%%%%%%%%%%%%%%%%%%%%%%%%%%%%%%%%%%%%%%%%%%%%%%%%%%%%%%%%%%
%%%%%%%%%%%%%%%%%%%%%%%%%%%%%%%%%%%%%%%%%%%%%%%%%%%%%%%%%%%%%%%%%%%%%%%%%%%%%%%%%%%%%%%%%%%%%%%%%%%%%%%%%%%%%%%%%
%%%%%%%%%%%%%%%%%%%%%%%%%%%%%%%%%%%
\section{Application to the discrete case}\label{sec:approx_discret}
%%%%%%%%%%%%%%%%%%%%%%%%%%%%%%%%%%%

\setcounter{equation}{0}

%%%%%%%%%%%%%%%%%%%%%%%%%%%%%%%%%%%%%%%%%%%%%%%%%
%\subsection{Notation in the discrete setting of local score}\label{sec:intro_biol}
%%%%%%%%%%%%%%%%%%%%%%%%%%%%%%%%%%%%%%%%%%%%%%%%%

\noi Recall that the r.v. $S_n$ and the Lindley process $U_n$  are associated with the sequence $(\epsilon_i)_{i\geq 1}$ via  \eqref{defsn} and \eqref{defuk} respectively. The process $(U_k)$ is a non negative Markov chain. In the case where $(\epsilon_i)_{i\geq 1}$ are symmetric Bernoulli r.v.'s (i.e. $P(\epsilon_i=\pm 1)=1/2$), then  $(U_k)$ takes its values in $\mathbb{N}$ and moves as a symmetric random walk in $\{1,2,\cdots\}$ and being at $0$, it either stays at this level with probability $1/2$ or jumps to $1$ with probability $1/2$.

\noi In general, the trajectory of $(U_k)$ can be decomposed in a succession of $0$ and excursions above $0$.  An excursion of $(U_k)$ starting at $g$ and ending at $d$ is a  process $(e(k),\; 0\leq k\leq \zeta)$, where
 $$e(0)=U(g)=0,\quad  e(\zeta)=U(d)=0,\quad \zeta:=d -g>0$$
 and $e(k):=U(g +k)>0, \ \mbox{{\rm for any } }  0< k <\zeta.$
 %
 %For simplicity we write $e$ instead of $e^U$.
  As mentioned in the Introduction, the local score $\overline{U}_n$ is the maximum of $(U_k)$ up to time $n$ and can be interpreted as  the maximum of all the heights of the   excursions up to time $n$. Namely 
\begin{equation}\label{scloc1}
    \overline{U}_n:=\max_{0\leq k\leq n}U_k, \quad n\geq 0.
\end{equation}

\noi  We are interested in the highest {\it complete} excursion up to time $n$. We proceed as in the continuous time setting introducing
 %%-------------------------------------------
\begin{equation}\label{dIn0}
\begin{array}{ccl}
g_n &:=&\max\left\lbrace k\leq n\ ;\ U_k=0\right\rbrace, \qquad  \quad U^*_n:=\dis \overline{U}_{g_n} =\max_{0\leq k\leq g_n}U_k,\\
 \\
f_n^* &:=&\dis \max\left\lbrace k\leq g_n ;\ U_k=U_n^*\right\rbrace,  \qquad g_n^*:=g_{f_n^*}=\max\left\lbrace k\leq f_n^*\ ;\ U_k=0\right\rbrace,\\
\\
d_n^* &:=&\inf\left\lbrace k\geq f_n^*\ ;\ U_k=0\right\rbrace,  \qquad \quad \theta_n^* :=  f^*_n-g^*_n.
\end{array}
\end{equation}

%%--------------and in view of notation of Section %\ref{sec:thal_results}, let $g_n$ (respectively $d_n$) be the left (resp. right) end point of the excursion associated %with %%$\overline{U}_n$. If $d_n\leq n$, then the excursion related to $\overline{U}_n$ is complete. In the case $d_n>n$, the highest complete excursion is the one associated with the maximal %%height $\overline{U}_{g_n}$.

%\noi Denote $e^*_n$ the "highest" complete excursion up to time $n$ and $g^* _n$ respectively $d^*_n$ its left (resp. %right) end point. The couple $(U^*_n,\theta^*_n)$ is then %%defined by
%
%%\begin{equation}\label{durho1}
%    U^*_n:=\max_{0\leq k\leq \zeta^*_n}e^*_n(k),\quad \theta^*_n:=\max\{0\leq k\leq \zeta^*_n ,\; e^*_n(k)=U^*_n\}
%\end{equation}
%
%where $\zeta^*_n:=d^*_n-g^*_n$. Then the r.v.s can be rewritten as

%%-------------------------------------------
%\begin{equation}\label{dIn0}
%\begin{array}{lcl}
%g_n =\max\left\lbrace k\leq n\ ;\ U_k=0\right\rbrace\\
%U^*_n=\overline{U}_{g_n} \\
%f_n^*:=\max\left\lbrace k\leq g_n ;\ U_k=U_n^*\right\rbrace\\
%g_n^*=g_{f_n^*}=\max\left\lbrace k\leq f_n^*\ ;\ U_k=0\right\rbrace\\
%d_n^* =\inf\left\lbrace k\geq g_n^*\ ;\ U_k=0\right\rbrace\\
%\theta_n^* =  f^*_n-g^*_n.
%\end{array}
%\end{equation}
%%--------------

\begin{figure}[h]\label{LabelNotationsUk}
\centering
\includegraphics[scale=0.3]{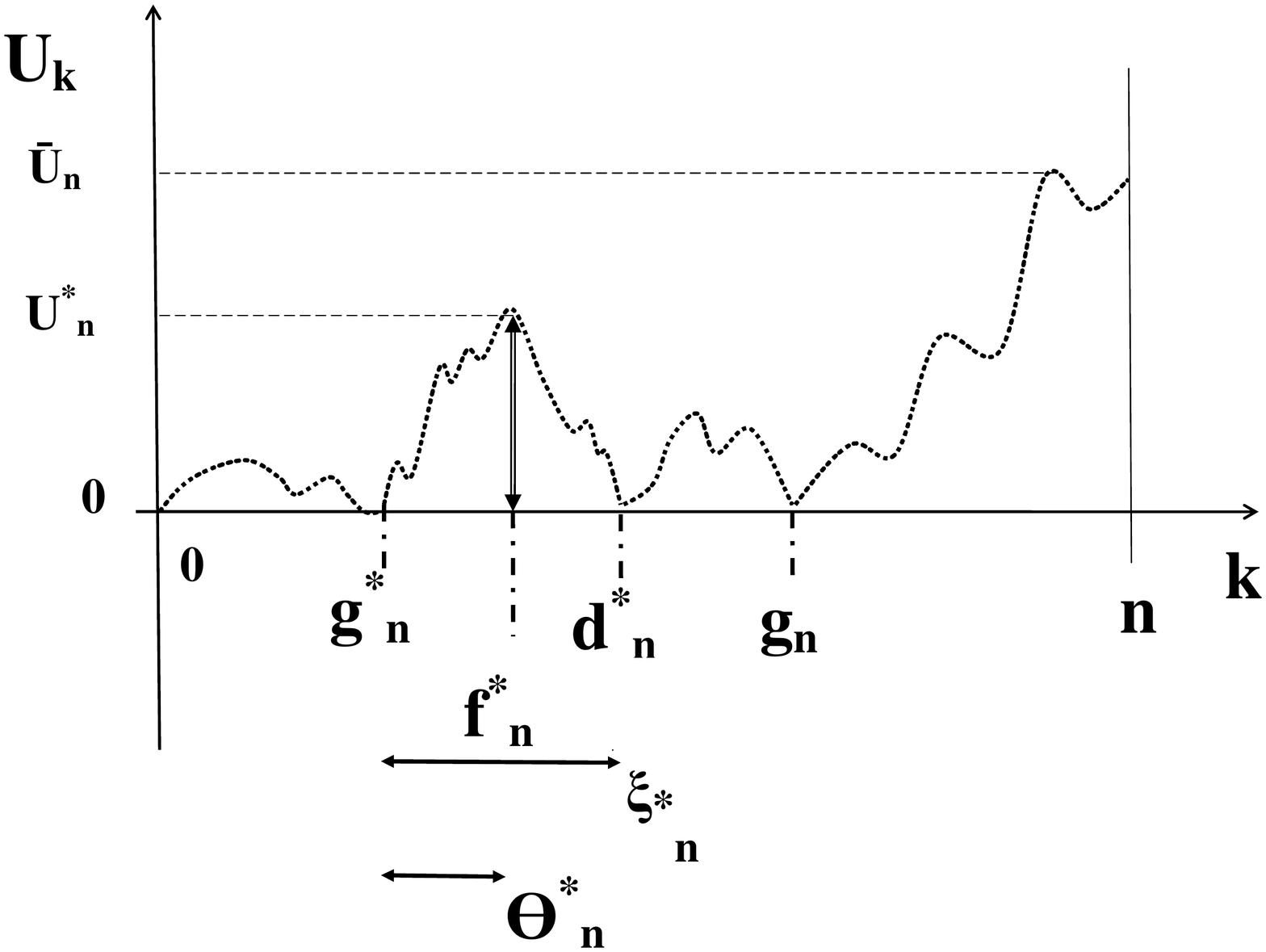}
\caption{Notation $U_n$}
\end{figure}
%%------------------------

\noi
In Section \ref{sec:link}, we define a continuous process $(U^{M}(t),t\geq 0)$ as the classical linear interpolation  of $(U_n,n\geq 0)$. We naturally introduce the highest high  $U^{M,*}(t)$ and length $\theta^{M,*}(t)$ of the complete excursion until time $t$ of $(U^{M}(t),t\geq 0)$. We conclude linking $(U^{M,*}(1),\theta^{M,*}(1))$ to $(U^*_M,\theta^*_M)$. Then we prove in Section \ref{sec:cv_proc} the convergence of $(U^{M,*}(t),\theta^{M,*}(t))$ to $(U^*(t),\theta^*(t))$. Since the distribution of $(U^*(t),\theta^*(t))$ has been computed in Section \ref{sec:thal_results}, we then get an approximation of the distribution of $(U^*_M,\theta^*_M)$.

%%%%%%%%%%%%%%%%%%%%%%%%%%%%%%%%%%%%%%%%%%%%%%%%%
\subsection{The linear interpolation of $(U_k)_k$}\label{sec:link}
%%%%%%%%%%%%%%%%%%%%%%%%%%%%%%%%%%%%%%%%%%%%%%%%%

%%A new continuous process $(U^{M}(t),t\geq 0)$
\noi We keep notation given  above and the one introduced in Section \ref{sec:thal_results}.
%subsection \ref{sec:intro_biol}.
Recall in particular that  $(B(t), \;t\geq 0)$ stands for a standard Brownian motion started at $0$ and  $(U(t), \;t\geq 0)$  is the reflected Brownian motion defined by \eqref{tr1}.
%
%
%\begin{equation}\label{trL1}
%    (U(t), \;t\geq 0) \ \overset{(d)}{= } \ \big(B(t)-\inf_{0\leq s\leq t}B(s),\; t\geq 0\big)  \ \overset{(d)}{= } \ %(\vert B(t)\vert , \;t\geq 0),
%\end{equation}
%
Let $M>0$ be a scale parameter which allows to obtain the convergence of the normalized random walk to the Brownian motion $(B(t))$ as $M\rightarrow \infty$ (see Section \ref{sec:cv_proc}). \noi
The classical continuous process
$\big(B^M(t),\; t\geq 0\big)$ associated with $(S_n)$ and normalizing factor $M$ is classically defined as $B^M\left(\frac{k}{M}\right)=\frac{1}{\sqrt M} S_k$ and for any $k$ such that 
$\frac{k}{M}\leq t\leq \frac{k+1}{M}$
$$B^M(t)=B^M\left(\frac{k}{M}\right)+M\left(t-\frac{k}{M}\right)\left(B^M\left(\frac{k+1}{M}\right)-B^M\left(\frac{k}{M}\right)\right).$$

\noi
We are interested here by the process $(U^M(t),\,t\geq 0)$
\begin{equation}\label{dUM}
    U^M(t)=B^M(t)-\min\limits_{s\leq t}B^M(s),\quad t\geq 0.
\end{equation}
Note that
%one has the natural identity in law
%\begin{equation}\label{trL2}
%   \big(|B^M(t)|, \;t\geq 0\big) \ \overset{(d)}{= } \ \big(B^M(t)-\min_{0\leq u\leq t}B^M(u),\; t\geq 0\big).
%\end{equation}
%
\begin{equation}\label{Ab1}
    U^M\left(\frac{k}{M}\right)=\frac{1}{\sqrt M} U_k,\quad k\geq 0
\end{equation}
where $(U_k)$ is the Lindley process associated with $(S_k)$ via \eqref{defuk}.\\

\noi
%In view of Corollary \ref{coUM} and Theorem \ref{T1}, define
We define the analog of r.v.s introduced in the discrete setting of Lindley process, see \eqref{dIn0} and \eqref{scloc1} in the continuous time setting of $\big(B^M(t)\big)$ 

%%-------------------------------------------
\begin{equation}\label{dIn1}
\begin{array}{ccl}
\overline{U}^M(t)&:=&\dis \sup_{0\leq s\leq t}U^M(s),  \qquad g^M(t) :=\sup\left\lbrace s\leq t\ ;\ U^M(s)=0\right\rbrace,\\
\\
U^{M,*}(t)&:=&\dis \overline{U}^M(g^M(t))= \sup_{0\leq s\leq g^M(t)}U^M(s),\\
\\
f^{M,*}(t)&:=&\dis \sup\left\lbrace r\leq g^M(t) ;\ U^M(r)=U^{M,*}(t)\right\rbrace, \\
\\
g^{M,*}(t)&:=&g^M(f^{M,*}(t))=\sup\left\lbrace r\leq f^{M,*}(t)\ ;\ U^M(r)=0\right\rbrace,\\
\\
d^{M,*}(t)&:=&\inf\left\lbrace s\geq f^{M,*}(t)\ ;\ U(s)=0\right\rbrace,  \qquad 
\theta^{M,*}(t):= f^{M,*}(t)-g^{M,*}(t).
\end{array}
\end{equation}
%%--------------

Using the definition \eqref{dIn0} of $\theta^*_M$ and $U^*_M$ we deduce easily that these r.v.'s can be expressed in terms of their analog in continuous time.

%%------------------
\begin{propo}\label{prop:link}
We have the following scaling properties
\begin{equation}\label{tr7}
  \frac{ \theta^*_M}{M}=\theta^{M,*}(1) \quad \trm{and} \quad \frac{U^*_M}{\sqrt{M }}=U^{M,*}(1).
\end{equation}
%%-------------
\end{propo}
%%-----------------------------

%%$$\frac{g_M}{M}=g^M(1) \quad \trm{and} \quad \frac{\overline{U}_M}{\sqrt{M }}=\overline{U}^M(1) \quad \trm{then} %%\quad  \frac{g^*_M}{M}=g^{M,*}(1) \quad \trm{and} \quad \frac{f^*_M}{M}=f^{M,*}(1) $$

%%%%%%%%%%%%%%%%%%%%%%%%%%%%%%%%%%%%%%%%%%%%%%%%%
%\subsection{Theoretical results of stochastic processes leading to Proposition \ref{P1}}\label{sec:tech_proc}
%%%%%%%%%%%%%%%%%%%%%%%%%%%%%%%%%%%%%%%%%%%%%%%%%

%%%%%%%%%%%%%%%%%%%%%%%%%%%%%%%%%%%%%%%%%%%%%%%%%
\subsection{Convergence of $\big(U^{M,*}(t),\theta^{M,*}(t)\big)$ to $\big(U^*(t),\theta^*(t)\big)$}\label{sec:cv_proc}
%%%%%%%%%%%%%%%%%%%%%%%%%%%%%%%%%%%%%%%%%%%%%%%%%

\noi The key ingredient of our convergence results is the Donsker Theorem, see Section 2.10 in \cite{Billingsley99}: the processus $\big(B^M(t), \; t\geq 0\big)$ converges weakly to the Brownian motion $\big(B(t),\; t\geq 0\big)$ when $M\to+\infty$. Using moreover \eqref{trL1} we get the following useful result.

%%---------------------------
\begin{propo}\label{coUM}
$\big(U^M(t), \; t\geq 0\big)$ converges weakly to $\big(U(t), \; t\geq 0\big)$.
\end{propo}
%%------------------------------------

\noi
Note that it is unclear that the map
 $\omega \mapsto \big(g^{M,*}(t),f^{M,*}(t), d^{M,*}(t), \theta^{M,*}(t)\big)$ defined from $U^*(t)$ is continuous.  Therefore the weak convergence of\\ 
$\big(g^{M,*}(t),f^{M,*}(t), d^{M,*}(t), \theta^{M,*}(t),U^{M,*}(t)\big)$ as $M\rightarrow \infty$ is not a straightforward consequence of Proposition \ref{coUM}.

% Note that it is unclear that the map
% $\omega \mapsto f^{M,*}(t)$ is continuous.  Therefore the weak convergence of
%$f^{M,*}(t)$, as $M\rightarrow \infty$, is not a straightforward consequence of Corollary \ref{coUM}.

%%%%%%%%%%%%%%%%%%%%%%%%%%%%%VERSION COURTE ???-----------------------------------------
%%-----------------------------
%\begin{theorem}\label{T1} Let $t>0$. Then, $\big(U^{M,*}(t),\; \theta^{M,*}(t)\big)$ converges weakly to %$\big(U^*(t),\; \theta^*(t)\big)$, as  $M \rightarrow \infty$.
%\end{theorem}
%%------------------------------------

%%-----------------------------
\begin{theorem}\label{T1} Let $t>0$.
%
%\begin{enumerate}
  %\item 
	The 5-uplet
$\big(g^{M,*}(t),f^{M,*}(t), d^{M,*}(t), \theta^{M,*}(t),U^{M,*}(t)\big)$ converges weakly to
$\big(g^*(t), f^*(t), d^*(t), \theta^*(t),  U^*(t)\big)$ as  $M\to\infty$ where the r.v.s $g^*(t)$, $f^*(t)$, $d^*(t)$, $\theta^*(t)$, $U^*(t))$ have been defined by relations \eqref{formule_U}-\eqref{tr8}.
  %\item In particular, $\big(U^{M,*}(t),\; \theta^{M,*}(t)\big)$ converges weakly to $\big(U^*(t),\; \theta^*(t)\big)$, as  $M \rightarrow \infty$.
%\end{enumerate}
 %
\end{theorem}

%\noi The proof of Theorem \ref{T1}  is postponed in Section \ref{app:tech_proc}.

%%%%%%%%%%%%%%%%%%%%%----------------------------------------------

%%-----------------------------
%\begin{propo}\label{P1} The couple
%$\big(U^{M,*}(t),\; \theta^{M,*}(t)\big)$ converges weakly to $\big(U^*(t),\; \theta^*(t)\big)$, as  $M %\rightarrow \infty$.
%\end{propo}
%%----------------------------------

%%%%%%%%%%%%%%%%%%%%%%%%%%%%%%%%%%%%%%%%%%%%%%%%%%%%%%%%%%%%%%%%%%%%%%%%%%%%%%%%%%%%%%%%%%%%%%%%%%%%
%%%%%%%%%%%%%%%%%%%%%%%%%%%%%%%%%%%%%%%%%%%%%%%%%%%%%%%%%%%%%%%%%%%%%%%%%%%%%%%%%%%%%%%%%%%%%%%%%%%%
%%%%%%%%%%%%%%%%%%%%%%%%%%%%%%%%%%%%%%%%%%%%%%%%%%%%%%%%%%%%%%%%%%%%%%%%%%%%%%%%%%%%%%%%%%%%%%%%%%%%

%\appendix
%%%%%%%%%%%%%%%%%%%%%%%%%%%%%%%
\section{Proofs}\label{Pro1}
%%%%%%%%%%%%%%%%%%%%%%%%%%%%%%%

\setcounter{equation}{0}

%\noi
%The theory of excursions of the Brownian motion and Propositions \ref{prop:prop1} and \ref{prop:prop2} are the main ingredients to prove Theorem \ref{th:th1}. 
%%Recall that $\big(L(t),\; t\geq 0\big)$ is the local time process at $0$ associated with $(B(t),\; t\geq 0)$, $\big(\tau_t,\; t\geq 0\big)$ is its right inverse and the r.v.s $U^*(t)$ and $\theta^*(t)$ have been defined by \eqref{formule_U} and \eqref{tr8}.

We follow the notation introduced in Sections \ref{sec:thal_results} and \ref{sec:approx_discret}.

%%%%%%%%%%%%%%%%%%%%%%%%%%%%%%%
\subsection{Proof of Proposition \ref{prop:prop1} }\label{app:th1}
%%%%%%%%%%%%%%%%%%%%%%%%%%%%%%%

We first link the distribution of $(U^*(t),\theta^*(t))$ to that of $(\overline{U}(\tau_s),\theta^*(\tau_s),\tau_s)$.

%%%%%%%%%----------------------
\begin{lemma}\label{lem:lem_f1}
 Let $f:\mathbb{R}_+^2\to\mathbb{R}$ be a bounded Borel function. Then
\begin{equation}\label{eq:lem_f1}
\mathbb{E}\big[f(U^*(t),\theta^*(t))\big]=\sqrt{\frac{2}{\pi}}\int_0^{+\infty}\mathbb{E}\left[f\big(\overline{U}(\tau_s),\theta^*(\tau_s)\big)
\frac{1}{\sqrt{t-\tau_s}} \1_{\{\tau_s<t\}}\right]\,ds.
\end{equation}
\end{lemma}
%%%%%%%%------------------------------

%%--------------------------------
\tbf{Proof} \
The real number $s=L(t)$ is the unique $s$ such that $\tau_{s_-}<t<\tau_s$.  Thus,
$$f\big(U^*(t),\theta^*(t)\big)=\sum\limits_{s\geq 0}\1_{\{\tau_{s_-}<t<\tau_s\}}f\big(\overline{U}(\tau_{s_-}),\theta^*(\tau_{s_-})\big)$$
since  $B(\tau_{s_-})=0$ implies that  $U^*(t)=\overline{U}(\tau_{s_-})$ and $\theta^*(t)=\theta^*(\tau_{s_-})$.

\noi Denote  $e_s$ the Brownian excursion
$$\begin{array}{lcl}
e_s(v):&=&\begin{cases}
B(\tau_{s_-}+v),\ 0\leq v\leq \tau_s-\tau_{s_-}\mbox{ for }\tau_s-\tau_{s_-}>0\\
[\delta] \ \mbox{ otherwise}
\end{cases}
\end{array}$$
and $\zeta(e_s):=\tau_s-\tau_{s_-}$ its lifetime. Since $\tau_s=\tau_{s_-}+\zeta(e_s)$,
$$\mathbb{E}\big[f(U^*(t),\theta^*(t))\big]=\mathbb{E}\left[\sum\limits_{s\geq 0} f\big(\overline{U}(\tau_{s-}),\theta^*(\tau_{s-})\big)\1_{\{\tau_{s-}<t<\tau_{s-}+\zeta(e_s)\}}\right].$$

\noi Applying Proposition 2.6 in \cite{RY99} (consequence of the Master Formula stated in Proposition 1.10, Chapter XII), one gets
$$\mathbb{E}\big[f(U^*(t),\theta^*(t))\big]=\mathbb{E}\left[\int_0^{+\infty}\left\lbrace\int f\big(\overline{U}(\tau_s),\theta^*(\tau_s)\big)\1_{\{\tau_s<t<\tau_s+\zeta(w)\}} n(dw)\right\rbrace\,ds\right],$$

\noi $n(dw)$ being a $\sigma$-finite measure on the set of all positive excursions. %In particular, by It\^o's description of Brownian excursions (see
According to Proposition 2.8, Chapter XII in \cite{RY99},  $n\big(\zeta(\omega)>\eps\big)=\sqrt{\frac{2}{\pi\eps}}$. Identity \eqref{eq:lem_f1} then follows.
%\noi from which we deduce
%%
%$$\mathbb{E}\big[f(U^*(t),\theta^*(t))\big]=\int_0^{+\infty}\mathbb{E}\left[f\big(\overline{U}(\tau_s),\theta^*(\tau_s)\big)
%\1_{\{\tau_s<t\}}\times\sqrt{\frac{2}{\pi}}\frac{1}{\sqrt{t-\tau_s}}\right]\,ds.$$
%%
\hfill $\blacksquare$

%\noi Using the scaling property of the process $U$ with index $1/2$, we state immediately the following lemma that reduces the problem to the particular case $s=1$.\\
Since for any $a>0$, the process $\left(U(sa)/\sqrt a; s\geq 0\right)$ is distributed as $\left(U(s); s\geq 0\right)$, we deduce the following scaling property

%%%%-------------
%\begin{lemma}\label{lem:lem_scale}
%Let $s>0$. Then
\begin{equation}\label{eq:scaling}
(\overline{U}(\tau_s),\theta^*(\tau_s),\tau_s)\overset{(d)}{=}(s\overline{U}(\tau_1),s^2\theta^*(\tau_1),s^2\tau_1), \quad s>0.
\end{equation}
%\end{lemma}
%%%---------------

%%%-----------
%\tbf{Proof} \
%According to  Corollary 2.2 of \cite{RY99}, $\tau_s$ is a stable subordinator of index $1/2$ and as a consequence $\tau_s\overset{(d)}{=}s^2\tau_1$.
%
%
%\noi Using  formula (10) in \cite{V91}, we have
%
%\begin{equation}\label{formule_scale_B}
%\left(\frac{1}{s}B(s^2r)\ ;\ 0\leq r\leq \frac{\tau_s}{s^2}\right)\overset{(d)}{=}\left(B(r)\ ;\ 0\leq r\leq \tau_1\right)
%\end{equation}
%
%\noi Since $B(\tau_s)=0$,
%\begin{eqnarray*}U^*(\tau_s)&=&\overline{U}(\tau_s)=\max\limits_{r\leq \tau_s}U(r)=\max\limits_{r\leq \tau_s}|B(r)|=\max\limits_{r'\leq \tau_s/s^2}|B(s^2 r')|=s\max\limits_{r'\leq \tau_s/s^2}\left(\frac{1}{s}B(s^2r')\right).\end{eqnarray*}
%
%\noi Similarly with $r=s^2r'$,
%\begin{eqnarray*}
%f^*(\tau_s)&=&\sup\{r<\tau_s\ ,\ |B(r)|=U^*(\tau_s)\}\\
%&=&\sup\{s^2r'<\tau_s\ ,\ |B(s^2r')|=U^*(\tau_s)\}\\
%&=&s^2\sup\left\lbrace r'<\frac{\tau_s}{s^2}\ ,\ \left|\frac{1}{s}B(s^2r')\right|=U^*(\tau_s)\right\rbrace
%\end{eqnarray*}
%
%\noi and
%\begin{eqnarray*}
%g^*(\tau_s)&=&\sup\{r<f^*(\tau_s)\ ,\ |B(r)|=0\}=s^2\sup\left\lbrace r'<\frac{f^*(\tau_s)}{s^2}\ ,\ \left(\frac{1}{s}B(s^2r')\right)=0\right\rbrace.
%\end{eqnarray*}
%
%\noi We conclude using $\theta^*(\tau_s)=f^*(\tau_s)-g^*(\tau_s)$.
%\hfill $\blacksquare$

\noi
Now we express the distribution of $\left(\overline{U}(\tau_1),\theta^*(\tau_1),\tau_1\right)$ in terms of that of $\left(\overline{U}(\tau_1),\tau_1\right)$.

\begin{lemma}\label{lem:lem_h1}
Let $h:\mathbb{R}_+^3\to\mathbb{R}$ be a bounded Borel function. Then
$$\mathbb{E}\left[h\left(\overline{U}(\tau_1),\theta^*(\tau_1),\tau_1\right)\right]=
\int_0^{+\infty}\mathbb{E}\left[h\big(x,x^2\xi,x^2(\xi+\xi')+\tau_1\big)\1_{\{\overline{U}(\tau_1)<x\}}\right]\,\frac{dx}{x^2}.$$
\end{lemma}

\tbf{Proof}\
It can be deduced from  Theorem 1 of \cite{V91} that
\begin{equation}\label{pt11}
  \mathbb{P}\big(\overline{U}(\tau_1)<x\big)=e^{-1/x} \quad  x>0.
\end{equation}
%
%See also Remark \ref{rPro1} below for a direct proof of \eqref{pt11}. 
Moreover conditionally on $\{\overline{U}(\tau_1)=x\}$,
\begin{enumerate}
  \item the r.v. $f^*(\tau_1)-g^*(\tau_1)$, $d^*(\tau_1)-f^*(\tau_1)$ and $\tau_1-d^*(\tau_1)+g^*(\tau_1)$ are independent;
  \item $f^*(\tau_1)-g^*(\tau_1)\overset{(d)}{=}d^*(\tau_1)-f^*(\tau_1)\overset{(d)}{=}T_x(R)$;
  \item $\tau_1-d^*(\tau_1)+g^*(\tau_1)$ is distributed as $\tau_1$ conditionally on $\{\overline{U}(\tau_1)<x\}$.
\end{enumerate}
%

%\noi
%We can see in Figure \ref{LabelNotationTau} the different variables of interest.
%
%
%%
%\begin{figure}[h]\label{LabelNotationTau}
%\centering
%\includegraphics[scale=0.3]{NotationsTau1.pdf}
%\caption{Notation with $\tau_1$}
%\end{figure}

\noi
%It is worth introducing the r.v.  $M_x$  having the same distribution as $\tau_1$ conditionally on $\{\overline{U}(\tau_1)<x\}$
%%
%$$\E\big(\phi(M_x)\big)=\E\big(\phi(\tau_1)\big|\overline{U}(\tau_1)<x\big)$$
%%
 %for any  bounded Borel function $\phi:\; [0,\infty[\rightarrow \mathbb{R}$. 
Now assume that $\widetilde{R}$ is distributed as  $R$ such that ($R$,$\widetilde{R}$) is independent of $U$. By the definition of $\theta^*(\tau_1)$ and a 
wise decomposition of $\tau_1$, we get
%\noi Since $\theta^*(\tau_1)=f^*(\tau_1)-g^*(\tau_1)$ and
%%
%$$\tau_1=\big(\tau_1-d^*(\tau_1)+g^*(\tau_1)\big)+\big(d^*(\tau_1)-f^*(\tau_1)\big)+\big(f^*(\tau_1)-g^*(\tau_1)\big),$$
%we deduce
\begin{eqnarray}
\mathbb{E}\big[h(\overline{U}(\tau_1),\theta^*(\tau_1),\tau_1)\big]&&= \int_0^{+\infty}\frac{e^{-1/x}}{x^2}\label{super_formule}\\
&&\mathbb{E}\left[h\big(x,T_x(R),T_x(R)+T_x(\widetilde{R})+\tau_1\big)| \overline{U}(\tau_1)<x\right]\,dx\nonumber.
\end{eqnarray}
%\begin{eqnarray}
%\mathbb{E}&&\big[h(\overline{U}(\tau_1),\theta^*(\tau_1),\tau_1)\big]=\int_0^{+\infty}\mathbb{E}\left[
%h\big(\overline{U}(\tau_1),\theta^*(\tau_1),\tau_1\big)\big|\overline{U}(\tau_1)=x)\right]f_{\overline{U}(\tau_1)}(x)dx\nonumber\\
%%&=&\int_0^{+\infty}\mathbb{E}\left[h\big(x,T_x(R),T_x(R)+T_x(\widetilde{R})+M_x\big)\right]f_{\overline{U}(\tau_1)}(x)dx.\nonumber\\
%&=&
 %\int_0^{+\infty}\mathbb{E}\left[h\big(x,T_x(R),T_x(R)+T_x(\widetilde{R})+\tau_1\big)| \overline{U}(\tau_1)<x\right]
 %\frac{e^{-1/x}}{x^2}\,dx\label{super_formule}.
%\end{eqnarray}
%
\noi The result is a direct consequence of \eqref{formule_xi} and  the scaling property
\begin{equation}\label{formule_T}
T_x(R)\overset{(d)}{=}x^2 T_1(R).
\end{equation}
%

%
%\begin{equation}
% \mathbb{E}\big[h(\overline{U}(\tau_1),\theta^*(\tau_1),\tau_1)\big]=
% \int_0^{+\infty}\mathbb{E}\left[h(x,T_x(R),T_x(R)+T_x(\widetilde{R})+\tau_1)| \overline{U}(\tau_1)<x\right]\frac{e^{-1/x}}{x^2}\,dx
%\end{equation}
%
%and $\mathbb{P}(U^*(\tau_1)<x)=\mathbb{P}(\overline{U}(\tau_1)<x)=e^{-1/x}$, we obtain the required result.
\hfill $\blacksquare$
%%%-------------------------

%\noi By Lemmas \ref{lem:lem_f1}, \ref{lem:lem_scale} and \ref{lem:lem_h1}, we now prove Proposition  \ref{prop:prop1}.

\tbf{Proof of Proposition \ref{prop:prop1}} Denote $\Delta:=\mathbb{E}\left[f(U^*(t),\theta^*(t))\right]$, where\\ $ f: \; [0,\infty[\times [0,\infty[ \rightarrow\mathbb{R}$ is a bounded Borel function.  According to  Lemma \ref{lem:lem_f1}, we have
$$\Delta=\sqrt{\frac{2}{\pi}}\int_{\R^2_+}f(y,y^2z)\mathbb{E}\left(\psi(y)\1_{\{y^2(\xi'+z)<t\}}\right)\frac{1}{y^2}p_{\xi}(z)\,dydz$$
where
\begin{eqnarray}\label{eq:def_psi}
\psi(y)&:=&\int_0^{+\infty}\frac{s\,ds}{\sqrt{t-y^2(z+\xi')-s^2\tau_1}}\1_{\{s<s_{\ast}\}}\\
&=&\frac{1}{\tau_1}\left[\sqrt{t-y^2(z+\xi')}-\sqrt{t-y^2(z+\xi')-s_\ast^2\tau_1}\right]
\end{eqnarray}
with $s_{\ast}=\frac{y}{\overline{U}(\tau_1)}\wedge\sqrt{\frac{t-y^2(z+\xi')}{\tau_1}}$ and $a\wedge b=\inf\{a,b\}.$

\noi %We explicit $\psi(y)$ as
It is easy to prove that on ${\left\lbrace s<s_{\ast}\right\rbrace}$, 
$$\psi(y)=\frac{1}{\tau_1}\left[\sqrt{\left(t-y^2(z+\xi')\right)_+}-\sqrt{\left(t-y^2\left\{z+\xi'-\frac{\tau_1}{\left(\overline{U}(\tau_1)\right)^2}\right\}\right)_+}\right].$$
Then \eqref{eq:dens_couple} follows.
\hfill $\blacksquare$
%%------------------------------

%%%--------------------
%\begin{rmk}\label{rem:rem25} Since $\tau_1\overset{(d)}{=}\frac{1}{G^2}$ where $G$ is a standard Gaussian r.v. (see e.g. \cite{RY99}),
%%
%\begin{equation}\label{loitau1}
    %\mathbb{P}(\tau_1\in ds)=\frac{1}{\sqrt{2\pi}}\frac{1}{s^{3/2}} e^{-\frac{1}{2s}}\1_{[0,+\infty[}(s)ds.
%\end{equation}
%%
%Consequently $\mathbb{E}\left(\frac{1}{\tau_1}\right)=+\infty$ and $E\Big(\frac{1}{\tau_1}\sqrt{t-y^2(\xi_1+\xi_2)} \big|\xi_1,\; \xi_2\Big)=\infty.$
%\end{rmk}
%%%-------------------------------

%%%%%%%%%%%%%%%%%%%%%%%%%%%%%%%%%
%%
\subsection{Proof of Proposition \ref{prop:prop2}}\label{sousp2}
%%
%%%%%%%%%%%%%%%%%%%%%%%%%%%%%%%%%%%%%%%%%%%%

Since the density function of $\overline{U}(\tau_1)$ is explicit (see \eqref{pt11}), that of
 $\big(\tau_1,\overline{U}(\tau_1)\big)$ will be determined once the conditional distribution of $\tau_1$ given $\overline{U}(\tau_1)$ will be known.  Our proof is based on the  study of the process $(\widehat{\lambda}(x),\; x>0)$ 
%which is supposed to be independent of $(U(t),\ t\geq 0)$ and 
such that conditionally on $\{\overline{U}(\tau_1)=x\}$,
\begin{equation}\label{formule_loi_tau_2}
\tau_1\overset{(d)}{=}\widehat{\lambda}(x).
\end{equation}
Obviously \eqref{formule_loi_tau_2} is equivalent to
\begin{equation}\label{Pro2}
    E\Big[f(\tau_1)g\big(\overline{U}(\tau_1)\big)\Big]=E\Big[f\Big(\widehat{\lambda}\big(\overline{U}(\tau_1)\big)\Big)g\big(\overline{U}(\tau_1)\big)\Big]
\end{equation}
for any bounded Borel functions $f$, $g: [0,\infty[\rightarrow\mathbb{R}$. We will show that $(\widehat{\lambda}(x),\; x>0)$  satisfies an equation which has a unique solution.

%%-----------------
\begin{lemma}\label{lem:lem_lambda}
Let $x>0$ and $n\geq 0$. Then,
\begin{equation}\label{formule_lambda_chap}
\widehat{\lambda}(x)\overset{(d)}{=}\Lambda_n+\widehat{\lambda}\left(\frac{1}{\frac{1}{x}+e_1+\cdots+e_{n+1}}\right)
\end{equation}
where
$$\Lambda_n:=x^2(\xi_1+\xi_2)+\sum\limits_{k=1}^n\frac{\xi_{2k+1}+\xi_{2k+2}}{\left(\frac{1}{x}+e_1+\cdots+e_k\right)^2}, \quad n\geq 0$$
with the classical convention $\sum\limits_1^0=0$.

\end{lemma}

\tbf{Proof} \  First we prove
$$\widehat{\lambda}(x)\overset{(d)}{=}x^2(\xi_1+\xi_2)+\widehat{\lambda}\left(\frac{1}{\frac{1}{x}+e_1}\right), \quad x>0.$$

\noi
Let $f_1$, $f_2: [0,+\infty[\to[0,+\infty[$ be two bounded Borel functions and
\begin{equation}\label{formule_A}
A:=\mathbb{E}\left[f_1(\tau_1)f_2\big(\overline{U}(\tau_1)\big)\right]=\int_0^{+\infty}e^{-1/x}f_2(x)\mathbb{E}[f_1(\widehat{\lambda}(x))]\,\frac{dx}{x^2}
\end{equation}
by \eqref{pt11} and (\ref{formule_loi_tau_2}). Applying formula (\ref{super_formule}) to $h(x_1,x_2,x_3)=f_1(x_3)f_2(x_1)$ leads to
\beq
A&=&\int_0^{+\infty}e^{-1/x}\mathbb{E}\left[f_1\big(T_x(R)+T_x(\widetilde{R})+\tau_1\big)\big|\overline{U}(\tau_1)<x\right] f_2(x)\,\frac{dx}{x^2}\\
%&=&\int_0^{+\infty}\frac{1}{x^2}\mathbb{E}\left[f_1\big(x^2(\xi_1+\xi_2)+\tau_1\big)\,\1_{\{\overline{U}(\tau_1)<x\}}\right]f_2(x)\,dx\\
&=&\int_0^{+\infty}\left(\int_0^x\frac{e^{-1/y}}{y^2}\mathbb{E}\left[f_1\left(x^2(\xi_1+\xi_2)+\widehat{\lambda}(y)\right)\right]\,dy\right) f_2(x)\,\frac{dx}{x^2}
\eeq
using (\ref{formule_T}), \eqref{pt11}  \eqref{formule_xi} and \eqref{Pro2}. Identifying with $(\ref{formule_A})$ implies
$$\mathbb{E}\left[f_1\big(\widehat{\lambda}(x)\big)\right]=e^{1/x}\int_0^x e^{-1/y}\mathbb{E}
\left[f_1\big(x^2(\xi_1+\xi_2)+\widehat{\lambda}(y)\big)\right]\,\frac{dy}{y^2}.$$

%\noi Let $Y$ the r.v. defined by $Y=\frac{1}{\frac{1}{x}+e_1}.$ Then $0<Y<x$ and for any $y\in]0,x[$, we get
%%
%$$P(Y<y)=P\left(\frac{1}{\frac{1}{x}+e_1}<y\right)=P\left(\frac{1}{x}+e_1>\frac{1}{y}\right)=
%P\left(e_1>\frac{1}{y}-\frac{1}{x}\right)=\exp\left(\frac{1}{x}-\frac{1}{y}\right).$$
%
%\noi Thus the density of $Y$ is $\frac{1}{y^2}e^{1/x}e^{-1/y}\,\1_{[0,x]}(y)$ and
%
\noi Let $Y$ be the r.v. defined by $Y=1\big/\left(\frac{1}{x}+e_1\right)$ whose density is obviously given by $e^{1/x-1/y}\,\1_{[0,x]}(y)/y^2$. Thus $\mathbb{E}\left[f_1\big(\widehat{\lambda}(x)\big)\right]$ can be rewritten as\\ $\mathbb{E}\left[f_1\big(x^2(\xi_1+\xi_2)+\widehat{\lambda}(Y)\big)\right]$ which means that 
$$\widehat{\lambda}(x) \overset{(d)}{=} x^2(\xi_1+\xi_2)+\widehat{\lambda}(Y)=x^2(\xi_1+\xi_2)+\widehat{\lambda}\left(\frac{1}{\frac{1}{x}+e_1}\right),\quad x> 0.$$
%
%\noi
%Now we iterate the procedure
%
%$$\widehat{\lambda}\left(\frac{1}{\frac{1}{x}+e_1}\right)\overset{(d)}{=}\left(\frac{1}{\frac{1}{x}+e_1}\right)^2(\xi_3+\xi_4)+\widehat{\lambda}\left(\frac{1}{\frac{1}{x}+e_1+e_2}\right)$$
%
%$$\widehat{\lambda}\left(\frac{1}{\frac{1}{x}+e_1+e_2}\right)\overset{(d)}{=}\left(\frac{1}{\frac{1}{x}+e_1+e_2}\right)^2(\xi_5+\xi_6)+\widehat{\lambda}\left(\frac{1}{\frac{1}{x}+e_1+e_2+e_3}\right)$$
%{\it etc}...\\
Iterating this procedure leads to \eqref{formule_lambda_chap}.
\hfill $\blacksquare$

\begin{lemma}\label{lem:lem_esp_lambda} For any $x>0$, $\mathbb{E}\left(\widehat{\lambda}(x)\right)=\frac 2 3 (x+x^2).$
\end{lemma}

\tbf{Proof} \ Using for instance Exercise (4.9) Chap VI in \cite{RY99} we get that
$$M(t):=\Big\{\cosh\big(\lambda|B(t)|\big)+b\sinh\big(\lambda|B(t)|\big)\Big\}\exp\left\lbrace-\frac{\lambda^2}{2}t-b\lambda L(t)\right\rbrace$$
is a local martingale for $\lambda>0$. Let $r>0$, $\dis b=-\frac{\cosh(\lambda r)}{\sinh(\lambda r)}$ and
$$\sigma_r:=\inf\{s\geq 0\ ;\ |B(s)|=r\}=\inf\{s>0\ ;\ U(s)=r\}.$$
The process $\left(M_{t\wedge \tau_1\wedge\sigma_r}; t\geq 0\right)$ being bounded, we can apply the stopping theorem to obtain
$\mathbb{E}(M(\tau_1\wedge\sigma_r))=\mathbb{E}(M(0)).$ It is clear that
$|B(\sigma_r)|=U(\sigma_r)=r,\quad B(\tau_1)=0$ and $L(\tau_1)=1.$
Our choice of  $b$ implies  that  $M(\sigma_r)=0$. Consequently,   $M(\tau_1\wedge\sigma_{r})=M(\tau_1)\,\1_{\{\tau_1<\sigma_r\}}$ and $e^{-b\lambda}\mathbb{E}\left[e^{-\lambda^2\tau_1/2}\,\1_{\{\tau_1<\sigma_r\}}\right]=1.$
Since $\{\tau_1<\sigma_r\}=\{\overline{U}(\tau_1)<r\}$, the previous identity can be rewritten as
\begin{equation}\label{Pro3}
  \mathbb{E}\left[e^{-\lambda^2\tau_1/2}\,\1_{\{\overline{U}(\tau_1)<r\}}\right]=e^{b\lambda}
  =\exp\left\{-\lambda\frac{\cosh(\lambda r)}{\sinh(\lambda r)}\right\},
\end{equation}
that leads to
\begin{equation}\label{Pro4}
  \mathbb{E}\left[e^{-\lambda\tau_1}\,\1_{\{\overline{U}(\tau_1)<r\}}\right]=\exp\left\lbrace-\sqrt{2\lambda}\;\frac{\cosh(r\sqrt{2\lambda})}
{\sinh(r\sqrt{2\lambda})}\right\rbrace=e^{-1/r}\left(1-\frac{2\lambda r}{3}+o(\lambda)\right).
\end{equation}
Taking the derivative at 0, we get
\begin{equation}\label{formule:esp_tau}
\mathbb{E}\left(\tau_1\1_{\{\overline{U}(\tau_1)<r\}}\right)=\frac{2r}{3}\; e^{-1/r}.
\end{equation}
%\begin{equation}\label{Pro4}
%  %\mathbb{E}\left[e^{-\lambda\tau_1}\,\1_{\{\overline{U}(\tau_1)<r\}}\right]=\exp\left\lbrace-\sqrt{2\lambda}\;\frac{\cosh(r\sqrt{2\lambda})}
%{\sinh(r\sqrt{2\lambda})}\right\rbrace.
%\end{equation}
%%
%Using classical analysis we get
%%
%$$\frac{u\cosh u}{\sinh  u}=\frac{u\left(1+\frac{u^2}{2}+o(u^2)\right)}{u\left(1+\frac{u^2}{6}+o(u^2)\right)}=\left(1+\frac{u^2}{2}\right)\left(1-\frac{u^2}{6}\right)+o(u^2)
%=1+\frac{u^2}{3}+o(u^2)\quad (u\rightarrow 0).$$
%
%\noi With $u=r\sqrt{2\lambda}$ and taking the limit $\lambda \rightarrow 0$ we obtain
%%
%$$\sqrt{2\lambda}\; \frac{\cosh(r\sqrt{2\lambda})}{\sinh(r\sqrt{2\lambda})}=\frac{1}{r}\left(1+\frac{2\lambda r^2}{3}+o(\lambda)\right)$$
%%
%%%=\frac{1}{r}\left[(\sqrt{2\lambda}r)\frac{{\rm ch}(\sqrt{2\lambda} r)}{{\rm sh}(\sqrt{2\lambda}r)}\right]
%%
%$$\exp\left\lbrace-\sqrt{2\lambda}\;\frac{\cosh(r\sqrt{2\lambda})}{\sinh(r\sqrt{2\lambda} )}\right\rbrace=e^{-1/r}\exp\left\lbrace-\frac{2\lambda r}{3}+o(\lambda)\right\rbrace=e^{-1/r}\left(1-\frac{2\lambda r}{3}+o(\lambda)\right).$$
%%
%Using \eqref{Pro4} we deduce
%%
%\begin{equation}\label{formule:esp_tau}
%\mathbb{E}\left(\tau_1\1_{\{\overline{U}(\tau_1)<r\}}\right)=\frac{2r}{3}\; e^{-1/r}.
%\end{equation}
%
\noi Let $\varphi$ be the function defined by  $\varphi(x):=\mathbb{E}\left[\widehat{\lambda}(x)\right].$
%=\mathbb{E}\big(\tau_1\big|\overline{U}(\tau_1)=x\big), \quad x\geq 0.$
%
Therefore, taking the conditional expectation with respect to $\overline{U}(\tau_1)$ in \eqref{formule:esp_tau} and using \eqref{pt11}, we have
%%
%\begin{eqnarray*}
%\mathbb{E}\big(\tau_1 \1_{\{\overline{U}(\tau_1)<r\}}\big)&=&\mathbb{E}\left(\mathbb{E}\big(\tau_1|\overline{U}(\tau_1)\big)\;\1_{\{\overline{U}(\tau_1)<r\}}\right)
%%=\int_0^r\frac{1}{x^2}e^{-1/x}\mathbb{E}\big(\widehat{\lambda}(x)\big)\,dx
%=\int_0^r\frac{1}{x^2}\; e^{-1/x}\varphi(x)\,dx.
%\end{eqnarray*}
%
%\noi Thus formula (\ref{formule:esp_tau}) can be rewritten as
%%
$\int_0^r e^{-1/x}\varphi(x)\,dx/x^2=\frac{2r}{3}\; e^{-1/r} \quad r>0,$ which conduces to\\ $\dis \varphi(x)=\frac{2}{3}\big(x+x^2\big)$. \hfill $\blacksquare$

%\frac{2x}{3}+\frac{2}{3}\,x^2
%$$\frac{1}{x^2}e^{-1/x}\varphi(x)=e^{-1/x}\left(\frac{1}{x^2}\frac{2x}{3}+\frac 2 3\right)$$
%after derivation and $\varphi(x)=\frac{2x}{3}+\frac{2}{3}\,x^2.$
%%-------------------------------

%%%--------------
%\begin{rmk}\label{rPro1} Note that if we take the limit $\lambda \rightarrow 0$ in \eqref{Pro4} we get
%%
%$$\mathbb{P}\big(\overline{U}(\tau_1)<r\big)=e^{-1/r},\quad \forall\; r>0$$
%%
%that is a direct proof of \eqref{pt11}.
%\end{rmk}
%%%-------------------------------------------

%%----------------------
\begin{lemma}\label{lem:lem_cv_lambda}
\begin{enumerate}
\item $\Lambda_n$ converges a.s. and in $L^1$ while $n\to+\infty$.
\item For any $x>0$, $\widehat{\lambda}\left(\frac{1}{\frac{1}{x}+e_1+\cdots+e_n}\right)$ converges to $0$ in $L^1$  while $n\to+\infty$.
\end{enumerate}
\end{lemma}
%%%%%%%%%%----------------------

%%%--------------------------------
\tbf{Proof}\
\ 1) Since all the r.v.s under concern are positive, $\Lambda_n$ converges a.s. while $n\rightarrow \infty$ to the positive r.v.
\begin{equation}\label{pl2x}
 \Lambda_\infty:=x^2(\xi_1+\xi_2)+\sum\limits_{k\geq 1}\frac{\xi_{2k+1}+\xi_{2k+2}}{\left(\frac{1}{x}+e_1+\cdots+e_k\right)^2}.
\end{equation}
One way to prove that $\Lambda_\infty$ is a.s. finite is to show that its expectation is finite.

\noi Note that (see \cite{BS02} p 463)
$$\mathbb{E}\left[e^{-\lambda\xi_1}\right]=\mathbb{E}\left[e^{-\lambda T_1(R)}\right]=\frac{\sqrt{2\lambda}}{\sinh(\sqrt{2\lambda})}=1-\frac{\lambda}{3}+o(\lambda),\quad  \lambda > 0.$$
%
%\noi An easy calculation leads to
%%
%$$\frac{u}{\sinh(u)}=\frac{u}{u+\frac{u^3}{3!}+o(u^3)}=1-\frac{u^2}{6}+o(u^2)\quad (u\rightarrow 0).$$
%%
\noi that conduces by derivation to $\mathbb{E}(\xi_1)=\frac 1 3.$ Now using the fact that $e_1+\cdots+e_k$ is $\gamma(k)$-distributed, we have successively
%
%$e_1+\cdots+e_k$ has a Gamma$(k)$ distribution,
%
\begin{eqnarray*}
\mathbb{E}(\Lambda_\infty)-\frac{2}{3}x^2&=&\mathbb{E}\left(\sum\limits_{k\geq 1}\frac{\xi_{2k+1}+\xi_{2k+2}}{\left(\frac{1}{x}+e_1+\cdots+e_k\right)^2}\right)\\ %=\sum\limits_{k\geq 1}\mathbb{E}\left(\frac{\xi_{2k+1}+\xi_{2k+2}}{\left(\frac{1}{x}+e_1+\cdots+e_k\right)^2}\right)\\
%&=&\sum\limits_{k\geq 1}\Big[\mathbb{E}(\xi_{2k+1})+\mathbb{E}(\xi_{2k+2})\Big]\;\mathbb{E}\left(\frac{1}{\left(\frac{1}{x}+e_1+\cdots+e_k\right)^2}\right)\\
&=&\frac 2 3 \sum\limits_{k\geq 1}\int_0^{+\infty}\frac{1}{\left(\frac{1}{x}+y\right)^2}\frac{y^{k-1}}{(k-1)!}\,e^{-y}\,dy
%&=&\frac 2 3 \int_0^{+\infty}\frac{1}{\left(\frac{1}{x}+y\right)^2}\left(\sum\limits_{k\geq 1}\frac{y^{k-1}}{(k-1)!}\right)\,e^{-y}\,dy\\
=\frac 2 3 \int_0^{+\infty}\,\frac{dy}{\left(\frac{1}{x}+y\right)^2}<+\infty
\end{eqnarray*}
which proves item 1 of Lemma \ref{lem:lem_cv_lambda}.

\vskip7pt \noi 2) Since $\widehat{\lambda}(y)\geq 0$, it is sufficient to check that
$\lim\limits_{n\to+\infty}\mathbb{E}\left[\widehat{\lambda}\left(\frac{1}{\frac{1}{x}+e_1+\cdots+e_n}\right)\right]=0.$

\noi As $0<\frac{1}{\frac{1}{x}+e_1+\cdots+e_n}\leq x$ and $\lim\limits_{n\to+\infty}\frac{1}{\frac{1}{x}+e_1+\cdots+e_n}=0\mbox{ a.s.}$ (by the Law of Large Numbers), the Lebesgue's dominated convergence theorem directly implies
 $$ \lim\limits_{n\to+\infty}\mathbb{E}\left(\frac{1}{\frac{1}{x}+e_1+\cdots+e_n}\right)=\lim\limits_{n\to+\infty}
 \mathbb{E}\left(\frac{1}{\left(\frac{1}{x}+e_1+\cdots+e_n\right)^2}\right)=0.$$
 %
%\noi The same arguments conduce to %$\lim\limits_{n\to+\infty}\mathbb{E}\left(\left(\frac{1}{\frac{1}{x}+e_1+\cdots+e_n}\right)^2\right)=0$.
\noi It remains to use Lemma \ref{lem:lem_esp_lambda} to get
$$\mathbb{E}\left[\widehat{\lambda}\left(\frac{1}{\frac{1}{x}+e_1+\cdots+e_n}\right)\right]=\frac 2 3 \left\lbrace\mathbb{E}\left(\frac{1}{\frac{1}{x}+e_1+\cdots+e_n}\right)+\mathbb{E}\left(\frac{1}{\left(\frac{1}{x}+e_1+\cdots+e_n\right)^2}\right)\right\rbrace.$$
and conclude the proof. \hfill $\blacksquare$
\subsection{Proof of Theorem \ref{th:th2}}\label{app:th2}
%%%%%%%%%%%%%%%%%%%%%%%%%%%%%%%

%\noi The goal is to give an explicit expression of $\Delta:=\mathbb{E}\left[f(U^*(t),\theta^*(t))\right]$ for a particular $f$. We express $\Delta$ in terms of the expectation of a  r.v. that can be simulated and derive his approximation using the  Monte Carlo algorithm. 
\noi We revisit the results of Sections \ref{app:th1} and \ref{sousp2}, keeping the notation introduced there. Interpreting the Lebesgue integral as an expectation in Lemma \ref{lem:lem_f1} gives: 
%that for any bounded Borel function $f:\mathbb{R}_+^2\to\mathbb{R}$,
%%%%%%%%%----------------------
%\begin{lemma}\label{lem:lem_f3}
 %Let $f:\mathbb{R}_+^2\to\mathbb{R}$ be a bounded Borel function. Then
\begin{equation*}\label{eq:lem_f3}
\mathbb{E}\big[f(U^*(t),\theta^*(t))\big]=\sqrt{\frac{\pi}{2}}\E\left[f\left(\sqrt{\frac{ t}{\tau_1}}\alpha_1 \overline{U}(\tau_1),\frac{ t}{\tau_1}\alpha_1^2 \theta^*(\tau_1)\right)\frac{1}{\sqrt{\tau_1}}\right].
\end{equation*}
%\end{lemma}
%%%%%%%%%----------------------

%\tbf{Proof} We apply Lemma \ref{lem:lem_f1}, the scaling property and the variable change $s=r\sqrt{\frac{t}{\tau_1}}$.
%\hfill $\blacksquare$

\noi
By the same reasoning, Lemma \ref{lem:lem_h1} can be modified as:
%%%%%%%%%----------------------
%\begin{lemma}\label{lem:lem_f4} 
%Let $h:\mathbb{R}_+^3\to\mathbb{R}$ be a bounded Borel function. Then
\begin{equation*}\label{eq:lem_f4}
\mathbb{E}\left[h\big(\overline{U}(\tau_1),\theta^*(\tau_1),\tau_1\big)\right]=\E\left[h\left(\frac{\overline{U}(\tau_1)}{\alpha_2},\frac{\overline{U}(\tau_1)^2}{\alpha_2^2}\xi,\frac{\overline{U}(\tau_1)^2}{\alpha_2^2}(\xi+\xi')+\tau_1\right)\frac{1}{\overline{U}(\tau_1)}\right].
\end{equation*}
%\end{lemma}
%%%%%%%%%----------------------
%\tbf{Proof} We apply Lemma \ref{lem:lem_h1} and the variable change $x=\frac{\overline{U}(\tau_1)}{y}$.
%\hfill $\blacksquare$
%
%
%
%\tbf{Proof of Theorem \ref{th:th2}}
Then the two previous equations and the following identity in law: $\overline{U}(\tau_1)\overset{(d)}{=}1/e_0'$ imply \eqref{15jui}.
%\beq
%\mathbb{E}\big[f(U^*(t),\theta^*(t))\big]&=&\sqrt{\frac{2}{\pi}}\mathbb{E}\left[f\big(\sqrt{\frac{t}{\tau_1}}\alpha_1 \overline{U}(\tau_1), \frac{t}{\tau_1}\alpha_1^2 \theta^*(\tau_1)\big)\frac{1}{\sqrt{\tau_1}}\right]\\
%\eeq
%We now state several identities in law. Using $(\overline{U}(\tau_1),\tau_1)\overset{(d)}{=}\left(1/e_0',\lambda(1/e_0')\right)$, we have
%
%\beq
%\frac{\overline{U}(\tau_1)}{\sqrt{\tau_1}}&\overset{(d)}{=}&\frac{\overline{U}(\tau_1)}{\sqrt{U(\tau_1)^2(\xi+\xi')+\tau_1\alpha_2^2}}
%\overset{(d)}{=}\frac{1}{\sqrt{\xi+\xi'+e_0'^2\alpha_2^2\lambda(1/e_0')}}\\
%\frac{\theta^*(\tau_1)}{\tau_1}&\overset{(d)}{=}&\frac{\overline{U}(\tau_1)^2\xi}{\overline{U}(\tau_1)^2(\xi+\xi')+\tau_1\alpha_2^2}
%\overset{(d)}{=}\frac{\xi}{\sqrt{\xi+\xi'+e_0'^2\alpha_2^2\lambda(1/e_0')}}\\
%\tau_1&\overset{(d)}{=}&\frac{\overline{U}(\tau_1)^2}{\alpha_2^2}(\xi+\xi')+\tau_1\overset{(d)}{=}\frac{1}{e_0^{'2}\alpha_2^2}(\xi+\xi')+\tau_1
%\overset{(d)}{=}\frac{1}{e_0^{'2}\alpha_2^2}\left(\xi+\xi'+e_0^{'2}\alpha_2^2\tau_1\right).
%\eeq
%Introducing $Z:=\xi+\xi'+e_0^{'2}\alpha_2^2\lambda(1/e_0')$ leads to the result.
\hfill $\blacksquare$

\subsection{Proof of Theorem \ref{th:marginal_theta}}\label{sec:th:marginal_theta}
%%%%%%%%%%%%%%%%%%%%%%%%%%%%%%%%%%%%%%%%%%%%%%%%%%%%%%%%%%%%%%%%%%%%%%

For any $a>0$, $b$ and $c\in \R$, we set
\begin{eqnarray*}
&&H(a,b):=\E\left [\frac{1}{\sqrt{b-a\xi}}\1_{\{b-a\xi>0\}}\right ]\\
&&\widehat{H}(a,b,c):=\E\left [\frac{1}{(a\xi-b)^{3/2}}\exp{\left(-\frac{c}{a\xi-b}\right)}\1_{\{a\xi-b>0\}}\right ].\\
\end{eqnarray*}

The proof of Theorem \ref{th:marginal_theta} is based on the following Lemma.

\begin{lemma}\label{lem:marginal_theta}
We have 
\begin{eqnarray}
&&H(a,b)=\frac{a}{b^{3/2}}\sum_{k\in\Z}|1+2k|\exp{\left (-\frac{(1+2k)^2}{2}\frac{a}{b}\right )}\label{H(a,b)}\\
&&\widehat{H}(a,b,c)
=\frac{\pi^{5/2}}{2a\sqrt{c}}\sum_{k\in\Z}(-1)^{k+1}k^2\exp{\left (-\frac{k^2\pi^2b}{2a}-|k|\pi\sqrt{\frac{2c}{a}}\right )}\label{H^(a,b,c)}.
\end{eqnarray}
\end{lemma}

\tbf{Proof} 1) Since $H(a,b)=0$ for $b\leq 0$, we assume from now on $b>0$. By \eqref{eq:dens_xi}
$$H(a,b)=\frac{1}{\sqrt{2\pi}}\sum_{k\in\Z}\left (-H_{1,k}+(1+2k)^2H_{2,k}\right )$$
with
\beq
H_{1,k}&:=&\int_0^{b/a}\frac{1}{u^{3/2}}\frac{1}{\sqrt{b-au}}\exp\left (-\frac{(1+2k)^2}{2u}\right )\,du,\\
H_{2,k}&:=&\int_0^{b/a}\frac{1}{u^{5/2}\sqrt{b-au}}\exp\left (-\frac{(1+2k)^2}{2u}\right )\,du.\\
\eeq
The change of variable $z=1/u-1/a$ in the above integrals gives
\beq
H_{1,k}&=&\sqrt{\frac{2\pi}{b}}\frac{1}{|1+2k|}\exp\left (-\frac{(1+2k)^2}{2}\frac{a}{b}\right),\\
H_{2,k}
&=&\left (\frac{a\sqrt{2\pi}}{b^{3/2}} \frac{1}{|1+2k|}+\frac{\sqrt{2\pi}}{\sqrt{b}}\frac{1}{|1+2k|^3}\right )\exp \left (-\frac{(1+2k)^2}{2}\frac{a}{b}\right ).
\eeq
From these relations we deduce the identity 
$-H_{1,k}+(1+2k)^2H_{2,k}=\frac{a\sqrt{2\pi}}{b^{3/2}}|1+2k|\exp\left ( -\frac{(1+2k)^2}{2}\frac{a}{b} \right )$ and finally \eqref{H(a,b)}.

2) Using \eqref{eq:dens_xi_2} and
%thus
%$$\widehat{H}(a,b,c)=\int_{b/a}^{\infty}\frac{\exp(-c/(au-b))}{(au-b)^{3/2}}\mu_2'(u)\,du.$$
an integration by parts lead to $\widehat{H}(a,b,c)=a\sum_{n\in\Z}(-1)^n\widehat{H}_n$
with
\begin{equation*}
\widehat{H}_n=\int_{b/a}^{\infty}\left (\frac{3}{2(au-b)^{5/2}}-\frac{c}{(au-b)^{7/2}}\right )\exp\left(-\frac{c}{au-b}-\frac{n^2\pi^2u}{2}\right)\,du.
\end{equation*}
With the change of variable $s=c/(au-b)$, we get
\beq
\widehat{H}_n
&=&\frac{1}{a}\left(\widehat{H}_n^1+\widehat{H}_n^2\right)\exp\left(-\frac{n^2\pi^2}{2}\frac{b}{a}\right)
\eeq
where
\beq
\widehat{H}_n^1&=&\frac{3}{2c^{3/2}}\int_0^{\infty}\sqrt{s}\exp\left(-s-\frac{n^2\pi^2c}{2a} \frac 1 s\right)\,ds,\\
\widehat{H}_n^2&=&-\frac{1}{c^{3/2}}\int_0^{\infty}s^{3/2}\exp\left(-s-\frac{n^2\pi^2c}{2a}\frac 1 s\right)\,ds.
\eeq
The Bessel functions $K_{\nu}$ admits the following integral representation (see formula (15) p 183 in\cite{Watson95}):
$$K_{\nu}(z)=\frac{1}{2}\left ( \frac{z}{2}\right )^{\nu}\int_{0}^{\infty}\frac{1}{s^{\nu+1}}\exp\left(-s-\frac{z^2}{4s}\right)\,ds.$$
Since $K_{-\nu}(z)=K_{\nu}(z)$ (see formula (8) p 79 in \cite{Watson95}, 
we obtain
$$\widehat{H}_n^1=\frac{3}{c^{3/2}}\left(\frac{z}{2}\right )^{3/2}K_{3/2}(z)\quad\mbox{and}\quad
\widehat{H}_n^2=-\frac{2}{c^{3/2}}\left(\frac{z}{2}\right )^{5/2}K_{5/2}(z)$$
where $z=|n|\pi\sqrt{\frac{2c}{a}}$. The functions $K_{3/2}$ and $K_{5/2}$ are explicit (see formula (12) p 80 in \cite{Watson95}): 
$$K_{3/2}(z)=\sqrt{\frac{\pi}{2z}}e^{-z}\left(1+\frac{1}{z}\right)
\quad\mbox{and}\quad
K_{5/2}(z)=\sqrt{\frac{\pi}{2z}}e^{-z}\left(1+\frac 3z+\frac{3}{z^2}.\right) $$
Then, we deduce $\widehat{H}_n^1+\widehat{H}_n^2=-\frac{\sqrt{\pi}}{4c^{3/2}}z^2e^{-z}$ and 
\eqref{H^(a,b,c)}.\hfill $\blacksquare$

\tbf{Proof of Theorem \ref{th:marginal_theta}} Using Proposition \ref{prop:prop1} and \eqref{H(a,b)}, for $0<x<t$, we get
%
%TRANSITION
%
\begin{eqnarray}\label{eqLem2DocPierre15avril2014}
f_{\theta^*(t)}(x)&=&\frac{1}{\sqrt{2\pi tx}}\int_{\R^3_+} \E\left [ H\left (\frac{x}{ty},1-\frac{x}{t}-\frac{u^2x}{ty}\lambda\left(\frac{1}{v}\right)\right )\1_{\left\{1-\frac{x}{t}-\frac{u^2x}{ty}\lambda\left(1/v\right)>0\right\}}\right ]\nonumber\\
&& \qquad \qquad \qquad \frac{u}{\sqrt{y}}e^{-v}\1_{\{u<v\}}p_{\xi}(y)\,du\, dv\, dy\nonumber\\
&=&\frac{1}{\sqrt{2\pi tx}}\int_{\R^2_+}ue^{-v}\1_{\{0<u<v\}}f_1(u,v)\,du\, dv
\end{eqnarray}
where 
$$f_1(u,v):=\int_0^{\infty}\E\left [ H\left ( \frac{x}{ty},1-\frac{x}{t}-\frac{u^2x}{ty}\lambda\left(\frac{1}{v}\right)\right )
 1_{\left\{1-\frac{x}{t}-\frac{u^2x}{ty}\lambda(\frac{1}{v})>0\right\}}\right ]\frac{p_{\xi}(y)}{\sqrt{y}}\,dy.$$
Using \eqref{H(a,b)} with $a=\frac{x}{ty}$ and $b=1-\frac{x}{t}-\frac{u^2x}{ty}\lambda\left(\frac{1}{v}\right)$, we get 
$$\frac{H()}{\sqrt{y}}=\frac{x\sqrt{t}}{\left [(t-x)y-u^2x\lambda(1/v)\right ]^{3/2}}\sum_{k\in \Z} |1+2k|
\exp{\left \{ - \frac{(1+2k)^2}{2}\frac{x}{(t-x)y-u^2x\lambda(1/v)}\right \}}.$$
Using the definition of the function $\widehat{H}$, we have:
$$f_1(u,v)=x\sqrt{t}\sum_{k\in\Z}|1+2k|\E\left [ \widehat{H}\left(t-x,u^2x\lambda(1/v),\frac{(1+2k)^2}{2}x\right)\right ].$$
Thus the density function of $\theta^*(t)$ can be written as follows, for $0<x<t$,
\begin{eqnarray}\label{eqLem4DocPierre15avril2014}
f_{\theta^*(t)}(x)&=&\frac{\sqrt{x}}{\sqrt{2\pi}}\int_{\R^2_+}\, ue^{-v}\1_{\{0<u<v\}}\\
&&\left \{ \sum_{k\in\Z}|1+2k|\E\left [\widehat{H}\left(t-x,u^2x\lambda(1/v),\frac{(1+2k)^2}{2}x\right)\right ]\right \}\,du\,dv.\nonumber
\end{eqnarray}

Set $a=t-x$, $b=u^2x\lambda(1/v)$ and $c=\frac{(1+2k)^2}{2}x$. From \eqref{H^(a,b,c)}, we get
\beq
&&\E\left[\widehat{H}\left(t-x,u^2x\lambda(1/v),\frac{(1+2k)^2}{2}x\right)\right]=\frac{\pi^{5/2}}{\sqrt{2x}(t-x)|1+2k|}\\
&&\sum_{n\in \Z}(-1)^{n+1}n^2\exp\left (-\pi |n|
 |1+2k|\sqrt{\frac{x}{t-x}}\right ) \E\left [ \exp{-\left (\frac{n^2\pi^2}{2(t-x)}u^2x\lambda(1/v)\right )}\right ] .
\eeq

Identity \eqref{Pro3} and item 2 of Proposition \ref{prop:prop2} give
$$\int_0^r\E \left(e^{-\mu\lambda(x)}\right)e^{-1/x}\,\frac{dx}{x^2}= \exp {\left \{-\sqrt{2\mu} \coth\left(r\sqrt{2\mu}\right)\right \}}$$
and a derivation with respect to $r$ leads to
\begin{equation}\label{eq:laplace}
\E [ \exp {\{-\mu\lambda (r)\}} ]=  r^2e^{1/r} \left ( \frac{2\mu}{\sinh^2 (r\sqrt{2\mu})}\right)\exp{\left \{ -\sqrt{2\mu}\coth(r\sqrt{2\mu})\right \}}.
\end{equation}
Taking $\mu=\frac{n^2\pi^2}{2(t-x)}u^2x$ we get
\begin{eqnarray}\label{eqLem6DocPierre15avril2014}
&&\E\left[\widehat{H}\left(t-x,u^2x\lambda(1/v),\frac{(1+2k)^2}{2}x\right)\right]=\frac{\pi^{9/2}u^2e^v\sqrt{x}}{\sqrt{2}|1+2k|(t-x)^2v^2}\\
&& \sum_{n\in\Z^*} (-1)^{n+1}n^4
\frac{
\exp{
\left \{
-\frac{\pi|n|\sqrt{x}}{\sqrt{t-x}}
\left ( 
|1+2k|+u \coth{\left(
\frac{\pi u|n|}{v}\sqrt{\frac{x}{t-x}}
\right)}
\right )
\right \}
}
}
{\sinh^2\left (\frac{\pi u|n|}{v}\sqrt{\frac{x}{t-x}}\right )}.\nonumber
\end{eqnarray}

Equations \eqref{eqLem4DocPierre15avril2014} and \eqref{eqLem6DocPierre15avril2014} imply that
$f_{\theta^*(t)}(x)=\sum_{k,n\in \Z^*}f_{k,n}(x)$
where 
\beq
f_{k,n}(x)&=&\sqrt{\frac{x}{2\pi}}|1+2k|\frac{\pi^{9/2}\sqrt{x}}{\sqrt{2}|1+2k|(t-x)^2}(-1)^{n+1}n^4\\
&&\int_{\R^2}\frac{u^3}{v^2}\frac{
\exp{
\left \{
-\frac{\pi|n|\sqrt{x}}{\sqrt{t-x}}
\left ( 
|1+2k|+u \coth{(
\frac{\pi u|n|}{v}\sqrt{\frac{x}{t-x}}
)}
\right )
\right \}
}
}
{\sinh^2\left (\frac{\pi u|n|}{v}\sqrt{\frac{x}{t-x}}\right )}\1_{\{u<v\}}\,du\,dv.
\eeq

Letting $u=vs$ ($v$ being fixed) and integrating with respect to $dv$, we obtain

\beq
f_{k,n}(x)&=&(-1)^{n+1}\pi |n| \sqrt{\frac{1}{x(t-x)}} \exp{\left (-\pi |n| |1+2k|\sqrt{\frac{x}{t-x}}\right )}\\
&&\times \int_0^1 \frac{\sinh\left (\pi s|n|\sqrt{\frac{x}{t-x}}\right )}{\cosh^2\left (\pi |n|\sqrt{\frac{x}{t-x}}\right )}\, ds\\
&=&(-1)^{n+1}\frac{1}{2x}\exp{\left (-\pi |n| |1+2k|\sqrt{\frac{x}{t-x}}\right )} \tanh^2\left (\pi |n|\sqrt{\frac{x}{t-x}}\right )
\eeq
and straightforward computation leads to
$$
\sum_{k\in \Z^*}f_{k,n}(x)= (-1)^{n+1}\frac{1}{2x}\tanh^2 \left (\pi |n|\sqrt{\frac{x}{t-x}}\right )\left ( \sinh\left (\pi |n|\sqrt{\frac{x}{t-x}}\right )\right )^{-1}$$
%with
%$$f_2(x)=\sum_{k\in\Z}\exp{\left (-\pi |n| |1+2k|\sqrt{\frac{x}{t-x}}\right )} =\left ( \sinh\left (\pi |n|\sqrt{\frac{x}{t-x}}\right )\right )^{-1}.$$
that finally conduces to \eqref{eq:marginal_theta}. \hfill $\blacksquare$

%%%%%%%%%%%%%%%%%%%%%%%%%%%%%%%%%%%%%%%%%%%%%%%%%%%%%%%%%%%%%%%%%%%%%%
\subsection{Proof of \eqref{eq:def_phi} in Theorem \ref{th:marginals}}\label{sec:th:marginals}
%%%%%%%%%%%%%%%%%%%%%%%%%%%%%%%%%%%%%%%%%%%%%%%%%%%%%%%%%%%%%%%%%%%%%%

%%------------
%\tbf{Proof of Theorem \ref{th:marginals}} 
It is clear that \eqref{eq:dist_U*_PY}-\eqref{eq:dist_b*} directly imply 
\beq
\P(U^*(t)>x)%&=&\P\left(\sqrt{t}\sqrt{g(1)}b^*>x\right)=\P\left(b^*>\frac{x}{\sqrt{tg(1)}}\right)\\
&=&\E\left(2\sum_{k\geq 1} (-1)^{k-1} e^{-\frac{2k^2x^2}{tg(1)}} \right)
%\eeq
%Since $g(1)<1$, we have $\frac{2k^2x^2}{tg(1)}>\frac{2x^2}{t}k^2$ that justifies the inversion between the expectation and the series. Hence by  \eqref{eq:dist_g1}
%\beq
%\P(U^*(t)>x)&=&2\sum_{k\geq 1} (-1)^{k-1} \E\left(e^{-\frac{2k^2x^2}{tg(1)}} \right)\\
=\frac{2}{\pi} \sum_{k\geq 1} (-1)^{k-1} \int_0^1 e^{-\frac{2k^2x^2}{ty}}\, \frac{dy}{\sqrt{y(1-y)}}.
\eeq
We take the $x$-derivative and we set $u=1/y-1$. \eqref{eq:def_phi} follows easily. \hfill $\blacksquare$
\subsection{Proof of Proposition \ref{prop:marginals}}
%%%%%%%%%%%%%%%%%%%%%%%%%%%%%%%%%%%%%%%%%%%%%%%%%%%%%%%%%%%%%%%%%%%%%%

\noi
Let us introduce $d_a:=\inf \{t\geq T_U(a),\ U(t)=0\}$. It is clear that $\left\{U^*(t)>a\right\}=\left\{d_a<t\right\}$. The process $\left(U(s+T_U(a))-a,\; 0\leq s\leq d_a-T_U(a)\right)$ is distributed as  
$\left(B(s),\; 0\leq s\leq T_B(-a)\right)$ and is independent of 
$\left(U(s),\; s\leq T_U(a)\right)$, then $d_a\overset{(d)}{=}T_U(a)+T_{\widehat{B}}(a)$. This shows \eqref{eq:inverse_2}. \hfill $\blacksquare$

%%%%%%%%%%%%%%%%%%%%%%%%%%%%%%%%%%%%%%%%%%%%%%%%%%%%%%%%%%%%%%%%%%%%%%
\subsection{Proof of Theorem \ref{T1}}\label{app:tech_proc}
%%%%%%%%%%%%%%%%%%%%%%%%%%%%%%%%%%%%%%%%%%%%%%%%%%%%%%%%%%%%%%%%%%%%%%

%\noi
%Recall that the $5$-uplet  $(g^{M,*}(t),d^{M,*}(t),f^{M,*}(t),\theta^{M,*}(t),U^{M,*}(t))$ has been defined by \eqref{dIn1}.
%Since $\theta^{M,*}(t)=f^{M,*}(t)-g^{M,*}(t)$, Theorem \ref{T1} will be proved as soon as the convergence of\\
 %$(g^{M,*}(t),d^{M,*}(t),f^{M,*}(t),U^{M,*}(t))$ holds.\\

%\noi
%Before proving the theorem, we return to the discrete setting.

\subsubsection{Auxiliary results in the discrete setting}\label{ssec:discrete}

Let us go back to the random walk defined by \eqref{defsn} and introduce for any integer $n_1>0$,
$$S'_k:=S_{n_1+k}-S_{n_1},\quad k\geq 0.$$

\begin{lemma}[Key Property]
\begin{enumerate}
\item $U_{n+1}=\max(U_n+\epsilon_{n+1},0)$.
\item Let $k$ be an integer such as $k>0$. Then
$$U_{n_1+i}>0\; \; \forall i\in\{0,\cdots,k\} \iff U_{n_1}>0\; \trm{and} \; U_{n_1}+S'_i>0 \; \; \forall i\in\{1,\cdots,k\}.$$
In such a case
$$U_{n_1+i}=U_{n_1}+S'_i\mbox{ for } 1\leq i\leq k.$$
\end{enumerate}
\end{lemma}

%\tbf{Proof} We have \begin{eqnarray*}
%U_{n+1 }&=&S_{n+1}-\min\limits_{0\leq i\leq n+1} S_i%=\max\limits_{0\leq i\leq n+1}(S_{n+1}-S_i)
%=\max(\max\limits_{0\leq i\leq n} (S_{n+1}-S_i),0)\\
%%&=&\max(\max\limits_{0\leq i\leq n} (S_n+\epsilon_{n+1}-S_i),0)
%&=&\max(\epsilon_{n+1}+\max\limits_{0\leq i\leq n}(S_n-S_i),0)
%=\max(\epsilon_{n+1}+U_n,0).
%\end{eqnarray*}
%Item (2) is a direct consequence of (1).
%%\item First assume $U_{n_1+i}>0\; \; \forall i\in\{0,\cdots,k\}$. By 1.,
%%
%%$\bullet$ $0<U_{n_1}$.
%%
%%$\bullet$ $0< U_{n_1+1}= U_{n_1}+\epsilon_{n_1+1}=U_{n_1}+S'_1$.
%%
%%$\bullet$ $\ldots$
%%%$\bullet$ $0< U_{n_1+2}= U_{n_1+1}+\epsilon_{n_1+2}=U_{n_1}+\epsilon_{n_1+1}+\epsilon_{n_1+2}=U_{n_1}+S'_2$ $\ldots$
%%
%%\noi\\
%%Conversely assume $U_{n_1}>0$ and $U_{n_1}+S'_i>0 \; \; \forall i\in\{1,\cdots,k\}.$ Still by 1.,
%%
%%$\bullet$ since $U_{n_1}+\epsilon_{n_1+1}>0$, $U_{n_1+1}= U_{n_1}+\epsilon_{n_1+1}>0$.
%%
%%$\bullet$ since $U_{n_1+1}+\epsilon_{n_1+2}>0$, $U_{n_1+2}= U_{n_1+1}+\epsilon_{n_1+2}>0$.
%%\end{enumerate}
%
%\hfill $\blacksquare$

\noi
Now consider
\begin{equation}\label{def:A1b}
N:=\{g^{*}_n<n_1,\ n_2<f^{*}_n<n_3,\ n_4<d^{*}_n<n_5,\ U^{*}_n\geq b\}
\end{equation}
where $0<n_1<\cdots<n_5<n$ are integers and see \eqref{dIn0} and \eqref{defuk} for the definition of the r.v.'s $g^{*}_n, f^{*}_n, d^{*}_n, U^{*}_n$ and $(U_k)$.\\
Define
$$\dis \overline{U}(m_1, m_2):=\max_{m_1\leq i\leq m_2}U_i,\; \; \;\dis \underline{U}(m_1,m_2):=\min_{m_1\leq i\leq m_2}U_i$$
and
$$n'_i:=n_i-n_1,\quad 2\leq i\leq 5, \quad n':=n-n_1.$$
\noi The event $D$ can be decomposed as
\begin{equation}\label{pTd2}
   N= N^1\cap N^2\cap N^3\cap N^4
\end{equation}
where
\begin{eqnarray}
% \nonumber to remove numbering (before each equation)
  N^1 &:=& \big\{U_k>0,\; n_1\leq k\leq n_4\big\}=\big\{\underline{U}(n_1,n_4)>0\big\}   \label{pTd3}\\
  N^2 &:=& \big\{\overline{U}(n_2, n_3)\geq \overline{U}(0, n_2)\vee b\big\} \label{pTd3b}\\
  N^3 &:=& \big\{ \overline{U}(n_2, n_3) > \overline{U}(n_3, n)\big\} \label{pTd4}\\
  N^4 &:=& \big\{\exists \; k,\;  U_k=0, \; n_4\leq k\leq n_5\big\}=\big\{\underline{U}(n_4,n_5)\leq 0\big\} \label{pTd5}
\end{eqnarray}
\noi
%In such a decomposition, $N$ depends on the set $\underline{U}(n_4,n_5)\leq 0$ and the r.v. $\underline{U}(n_4,n_5)$ has an atom at $0$ which will lead us to a problem in the sequel. This explains why we rewrite the set $D$ into another form.
Now note that
\begin{equation}\label{cond_pos}
U_{n_1}+S'_i>0 \iff S_{n_1+k}-\underline S(0,n_1)>0.
\end{equation}

\begin{figure}[h]\label{LabelSuiteTi}
\centering
\includegraphics[scale=0.3]{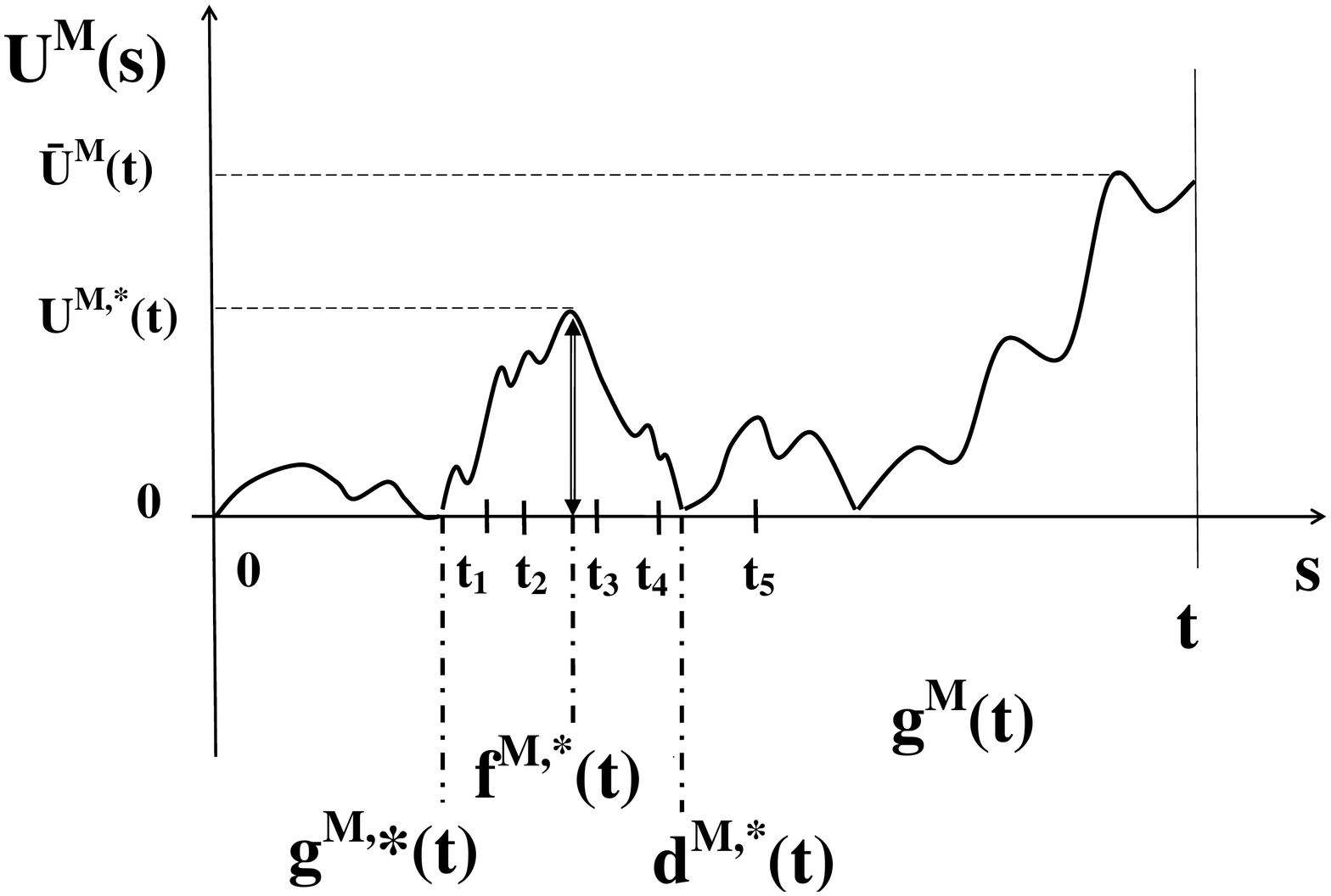}
\caption{Sequence of $t_i$}
\end{figure}

Moreover by the definitions of the $n_j$'s, one has $U_{n_j}>0$ and $S_{i+n_j} -\underline{S}(0,n_j)>0 \; \forall i=1\ldots n_5.$ Consequently, we successively have
\begin{eqnarray}
% \nonumber to remove numbering (before each equation)
  N ^1&=& \big\{\underline{S'}(0, n'_4)>-U_{n_1}\big\} \\
 % &&\nonumber\\
  N ^2&=& \Big\{\overline{S'}(n'_2, n'_3)\geq -U_{n_1}+\max\big[\overline{U}(0,n_1),\; b, \; \overline{S'}(0,n'_2)+U_{n_1}\big]\Big\} \label{pTd6}\\
  %&&\nonumber\\
  N ^3&=& \Big\{\overline{S'}(n'_2, n'_3)\geq -U_{n_1}+\overline{U}(n_3,n)\Big\} \label{pTd7}\\
 %&&\nonumber\\
  N ^4&=& \big\{\underline{S'}(n'_4, n'_5)\leq -U_{n_1}\big\}. \label{pTd8}
 \end{eqnarray}
\noi The above equalities  can  be directly read on Figure \ref{LabelSuiteTi} (a dash line representing a level that could not be crossed by the process).

\noi We want to express $N^{3}$ in terms of
$T'_k:=S_{t_3+k}-S_{t_3}.$
We have
$$\underline{S}(0,n_3+k)=\min\big\{ \underline{S}(0,n_3),\; S_{n_3}+\underline{T'}(0,k)\big\}$$
and
$$
\begin{array}{ccl}
U_{n_3+k}&=&S_{n_3+k}-\underline{S}(0,n_3+k)
= T'_{k}+\max \big\{U_{n_3},\; -\underline{T'}(0,k)\big\}.
\end{array}
$$
As a result
$$N^{3} = \Big\{ \overline{S'}(n_2',n_3')>-U_{n_1}+ \max_{0\leq k\leq n-n_3}\Big[
T'_k+\max \big\{U_{n_3},\;-\underline{T'}(0,k)\big\}\Big]\Big\}.$$

\subsubsection{Back to the continuous case}\label{ssec:cont}

\vskip5pt\noi\\ {\bf 1)} Let $t_1,t_2,\ldots,t_5$ be positive real numbers such that $0<t_1<\cdots<t_5$ and $b>0$.
Let us introduce
\begin{equation}\label{def:A1}
A_M^1=\{g^{M,*}(t)<t_1,\ t_2<f^{M,*}(t)<t_3,\ t_4<d^{M,*}(t)<t_5,\ U^{M,*}(t)>b\}
\end{equation}
where $g^{M,*}(t)$, $f^{M,*}(t)$, $d^{M,*}(t)$ and $U^{M,*}(t))$ have been defined by \eqref{dIn1}. The goal is to show
\begin{equation}\label{prT1a}
    \lim_{M\rightarrow \infty}P(A_M^1)=P(A^1)
\end{equation}
where
\begin{equation}\label{prT1c}
    A^1:=\{g^*(t)<t_1,\ t_2<f^*(t)<t_3,\ t_4<d^*(t) <t_5,\ U^*(t)>b\}
\end{equation}
and  the r.v.'s $U^*(t)$, $g^*(t), f^*(t)$ and $d^*(t)$ have been defined by \eqref{formule_U}-\eqref{formule_d}.

\vskip5pt\noi\\ {\bf 2)} In view of the discrete case, let us consider the sets of dyadic points
$$D=\underset{m\in \N}{\bigcup} D_m \quad \trm{where} \quad D_m=\left\{\frac{k}{2^m},\quad k\in\{0,1,\ldots\}\right\}.$$

\noi
Since $D$ is dense in $\R$ and $D_{n}\subset D_{m}$ as soon as $n\leq m$, we can choose without loss of generality positive integers $L_0$, $l$ and $l_i$ for $i=1\ldots 5$ such as
$$t_i=\frac{l_i}{2^{L_0}},\ 1\leq i\leq 5,\quad t=\frac{l}{2^{L_0}}.$$

Recall that $\big( U^M(t),\; t\geq 0\big)$ is the continuous process defined by \eqref{dUM} and the linear interpolation of $\dis \left(\frac{1}{\sqrt{M}}U_k,\; k\geq 0\right)$.

\vskip5pt\noi\\ {\bf 3)} For any  continuous function $\omega:\; [0,\infty[\rightarrow \mathbb{R}$, we denote
\begin{equation}\label{PrT1e}
\overline \omega(u, v):=\max_{u\leq r\leq v} \omega(r),\quad \underline{\omega}(u, v):=\min_{u\leq r\leq v}\omega(r),\quad 0\leq u \leq v.
\end{equation}

Following the procedure presented in the discrete case, the event $A_M^1$ can be decomposed as
\begin{equation}\label{pTd2b}
   A_M^1= A_M^{1,1}\cap A_M^{1,2}\cap A_M^{1,3}\cap A_M^{1,4}
\end{equation}
\noi\\
where for $i=1,\ldots,4$ $A_M^{1,i}$ is the analog of $N_i$ obtained by replacing $\underline{U}$ (resp. $\overline{U}$, $n_i$, $i=1,\ldots, 5$) by 
$\underline{U}^M$ (resp. $\overline{U}^M$, $t_i$, $i=1,\ldots, 5$).
%%
%\begin{eqnarray}
%% \nonumber to remove numbering (before each equation)
  %A_M^{1,1} &:=& \big\{U^M(t)>0,\; t_1\leq t\leq t_4\big\}=\big\{\underline{U}^M(t_1,t_4)>0\big\}  \label{pAd3}\\
  %A_M^{1,2} &:=& \big\{\overline{U}^M(t_2, t_3)\geq \overline{U}^M(0, t_2)\vee b\big\} \label{pAd3b}\\
  %A_M^{1,3} &:=& \big\{ \overline{U}^M(t_2, t_3) > \overline{U}^M(t_3, t)\big\} \label{pAd4}\\
  %A_M^{1,4} &:=& \big\{\exists \; s,\;  U^M(s)=0, \; t_4\leq t\leq t_5\big\}=\big\{\underline{U}^M(t_4,t_5)\leq 0\big\} .\label{pAd5}
%\end{eqnarray}

\noi\\
By Corollary \ref{coUM}, $\underline{U}^M(t_4,t_5)$ converges weakly to $\underline{U}(t_4,t_5)$, as $M\rightarrow \infty$. Thus we want to study the limit when $M$ goes to infinity and apply the following lemma

\begin{lemma}
Let $(\xi^M)$ be a sequence of r.v.'s valued in $\mathbb{R}^d$ and converging weakly to $\xi$ when $M\rightarrow \infty$. Then the Porte-Manteau's lemma (see e.g. \cite{Billingsley99}) asserts that for any Borel $\Lambda$ in $\mathbb{R}^d$,
\begin{equation}\label{pT1}
    \lim_{M\rightarrow \infty } \mathbb{P}(\xi^M\in \Lambda)=P(\xi\in \Lambda)
\end{equation}
if $\mathbb{P}(\xi \in \partial\Lambda)=0$.
\end{lemma}

Unfortunately, the distribution of $\underline{U}(t_4,t_5)$ (being bounded below by 0) has an atom at $0$; therefore we cannot conclude directly that $\dis \lim_{M\rightarrow \infty}\mathbb{P}\big(\underline{U}^M(t_4,t_5)= 0\big)=\mathbb{P}\big(\underline{U}(t_4,t_5)= 0\big)$. This is the reason why we will introduce the processes $W$ and $Z$ in the sequel.

\vskip5pt\noi\\ {\bf 4)} We follow now the procedure developed in section \ref{ssec:discrete}. It is worth introducing
$t_i':=t_i-t_1$, $i\in\{2, 3, 4, 5\}$, $t':=t-t_1$
and $W^M$ the process
$$W^M(s):=B^M(t_1+s)-B^M(t_1),\quad s\geq 0.$$
Note that the process $\big(W^M(s), s\geq 0\big)$ is the linear interpolation of $\dis \Big(\frac{1}{\sqrt{M}}(S_{k+n_1}-S_{n_1}, k\in \mathbb{N}\Big)$. We deduce from the previous step that $A_M^1=A_M^2$ where
$$A_M^2:=A_M^{2,1}\cap A_M^{2,2}\cap A_M^{2,3}\cap A_M^{2,4}$$
and
$$\begin{array}{ccl}
A_M^{2,1} &:=& \big\{\underline{W}^M(0,t_4')>-U^M(t_1)\big\}\\
%&&\\
A_M^{2,2} &:= & \Big\{ \overline{W}^M(t_2',t_3')\geq -U^M(t_1)+\max\big[\overline{U}^M(0,t_1),\; b,\; \overline{W}^M(0,t_2')+U^M(t_1)\big]\;\Big\} \\
%&&\\
%A_M^{2,3}  &:=& \Big\{ \overline{W}^M(t_2',t_3')> -U^M(t_1)+ \overline{U}^M(t_3, t)\Big\}\\
A_M^{2,3}   &:=& \Big\{ \overline{W}^M(t_2',t_3')>-U^M(t_1)+ \underset{0\leq u\leq t-t_3}{\max} \Big[
Z^M(u)+\max \big\{U^M(t_3),\;-\underline{Z}^M(0,u)\big\}\Big]\Big\}\\
%&&\\
A_M^{2,4}  &:=& \big\{\underline{W}^M(t_4',t_5')\leq -U^M(t_1)\big\}.
\end{array}
$$

%\noi We would like to express $A_M^{2,3}$ with
%%
%$$Z^M(u):=B^M(t_3+u)-B^M(t_3),\quad u\geq 0.$$
%%
%First, we have
%%
%$$\underline{B}^M(0,t_3+u)=\min\big\{ \underline{B}^M(0,t_3),\; B^M(t_3)+\underline{Z}^M(0,u)\big\}.$$
%%
%Second we deduce
%%
%$$
%\begin{array}{ccl}
%U^M(t_3+u)&=&B^M(t_3+u)-\underline{B}^M(0,t_3+u) \\
%&&\\
%& =& Z^M(u)+\max \big\{U^M(t_3),\; -\underline{Z}^M(0,u)\big\}.
%\end{array}
%$$
%%
%As a result
%%
%$$A_M^{2,3} = \Big\{ \overline{W}^M(t_2',t_3')>-U^M(t_1)+ \max_{0\leq u\leq t-t_3}\Big[
%Z^M(u)+\max \big\{U^M(t_3),\;-\underline{Z}^M(0,u)\big\}\Big]\Big\}$$

%

\vskip5pt\noi\\ {\bf 5)} To conclude the proof by taking the limit in $M$, it remains to express the limit subsets in the same way. In that view, let us introduce
$$W(s):=B(t_1+s)-B(t_1),\quad Z(s):=B(t_3+s)-B(t_3),\quad s\geq 0$$
and
$$A^2:=A^{2,1}\cap A^{2,2}\cap A^{2,3}\cap A^{2,4}$$
with
$$\begin{array}{ccl}
A^{2,1} &:=& \big\{\underline{W}(0,t'_4)>-U(t_1)\big\}\\
%&&\\
A^{2,2} &:= & \big\{ \overline{W}(t'_2,t'_4)\geq -U(t_1)+\max\big[ \overline{U}(0,t_1),\;b,\; \overline{W}(0, t'_2)+U(t_1)\big]\;\big\} \\
%&&\\
A^{2,3}  &:=& \dis \Big\{ \overline{W}(t_2',t_3')>-U(t_1)+ \max_{0\leq u\leq t-t_3}\Big[
Z(u)+\max \big\{U(t_3),\; -\underline{Z}(0,u)\big\}\Big]\Big\}\\%&&\\
A^{2,4}  &:=& \big\{\underline{W}(t_4',t_5')\leq -U(t_1)\big\}.
\end{array}
$$
Recall that for any $u>0$,  the r.v.s $\dis \max_{0\leq r\leq u}B(r)$ and $\dis \min_{0\leq r\leq u}B(r)$ have a density function. Therefore we can apply \eqref{pT1} to get
$$\lim_{M\rightarrow \infty}\mathbb{P}(A^{2,M})=\mathbb{P}(A^2).$$
As done in the discrete setting, we deduce that $A^2 =A^1$ where $A^1$ has been defined by \eqref{prT1c}. It is now clear that \eqref{prT1a} follows.

\hfill $\blacksquare$

%%%%---------------
%\begin{rmk}\label{rmpb1} According to Corollary \ref{coUM}, $\underline{U}^M(t_4,t_5)$ converges weakly to $\underline{U}(t_4,t_5)$, as $M\rightarrow \infty$. Unfortunately, the distribution of $\underline{U}(t_4,t_5)$ has an atom at $0$, therefore we cannot conclude directly that $\dis \lim_{M\rightarrow \infty}\mathbb{P}\big(\underline{U}^M(t_4,t_5)= 0\big)=\mathbb{P}\big(\underline{U}(t_4,t_5)= 0\big)$. This explains why we have introduced the processes $W$ and $Z$.
%\end{rmk}
%%%-----------

\bibliographystyle{plain}
\bibliography{biblio_sequence}

\end{document}